\documentclass[final,onefignum,onetabnum]{siamart190516}

\usepackage{lipsum}
\usepackage{amsfonts}
\usepackage{graphicx}
\usepackage{epstopdf}
\usepackage{algorithmic}
\usepackage{amsopn}
\DeclareMathOperator{\diag}{diag}

\usepackage{amsmath}
\usepackage{amssymb}
\usepackage{fancyhdr}
\usepackage{color}
\usepackage[caption=false]{subfig}
\usepackage{multirow}
\usepackage{rotating}
\usepackage{tabularx}
\usepackage{ltablex}
\usepackage{booktabs}
\usepackage{listings}
\usepackage{nicematrix}
 \allowdisplaybreaks[4]
\usepackage{longtable}
\newcommand{\tr}{^{\sf T}}
\newsiamremark{remark}{Remark}
\newsiamremark{hypothesis}{Hypothesis}
\crefname{hypothesis}{Hypothesis}{Hypotheses}
\newsiamthm{claim}{Claim}
\title{A mechanism of three-dimensional quadratic termination for the gradient method with applications\thanks{This work was supported by the National Natural Science Foundation of China (nos. 12021001, 11991021, 11991020, and 12071108) and the Strategic Priority Research Program of Chinese Academy of Sciences (no. XDA27000000).}}

\author{Ya-kui Huang\thanks{Institute of Mathematics, Hebei University of Technology, Tianjin 300401, China
  (\email{hyk@hebut.edu.cn}).}
\and Yu-Hong Dai\thanks{Corresponding author. LSEC, Academy of Mathematics and Systems Science, Chinese Academy of Sciences, Beijing 100190, China; School of Mathematical Sciences, University of Chinese Academy of Sciences, Beijing
	100049, China
  (\email{dyh@lsec.cc.ac.cn}, \url{http://lsec.cc.ac.cn/\string~dyh/}).}
\and Xin-Wei Liu\thanks{Institute of Mathematics, Hebei University of Technology, Tianjin 300401, China
  (\email{mathlxw@hebut.edu.cn}).}
}

\begin{document}

\maketitle

\begin{abstract}
Recent studies show that the two-dimensional quadratic termination property has great potential in improving performance of the gradient method. However, it is not clear whether  higher-dimensional quadratic termination leads further benefits.  In this paper, we provide an affirmative answer by introducing a mechanism of three-dimensional quadratic termination for the gradient method. A novel stepsize is derived from the mechanism such that a family of delayed gradient methods equipping  with the novel stepsize have the three-dimensional quadratic termination property. When applied to the Barzilai--Borwein (BB) method, the novel stepsize does not require the use of any exact line search or the Hessian, and can be computed by stepsizes and gradient norms in previous iterations.    Using long BB steps and some short steps associated with the novel stepsize in an adaptive manner, we develop an efficient gradient method for quadratic optimization and further extend it to general unconstrained optimization. Numerical experiments show that the three-dimensional quadratic termination property can significantly improve performance of the BB method, and the proposed method outperforms gradient methods that use stepsizes with the two-dimensional quadratic termination property.
\end{abstract}



\pagestyle{myheadings}
\thispagestyle{plain}
\markboth{Y.-K. HUANG, Y.-H. DAI and X.-W. LIU}
{3D QUADRATIC TERMINATION FOR THE GRADIENT METHOD}


\section{Introduction}
\label{intro}
As one of the most popular methods for solving smooth unconstrained optimization, the gradient method simply updates iterates by
\begin{equation}\label{eqitr}
	x_{k+1}=x_k-\alpha_kg_k,
\end{equation}
where $g_k=\nabla f(x_k)$ is the gradient of the objective function $f(x): \mathbb{R}^n \to \mathbb{R}$ at $x_k$ and $\alpha_k>0$ is the stepsize. A fundamental result of the gradient method \eqref{eqitr} is the quadratic termination property, i.e., when $f(x)$ is a strictly convex quadratic function, the exact minimizer can be achieved in finite iterations if the stepsizes sweep all inverse eigenvalues of the Hessian \cite{lai1981some}. However, such a quadratic termination property is generally not possible because computing all eigenvalues of a given $n\times n$ matrix requires $\mathcal{O}(n^3)$ floating point operations (flops) \cite{pan1999complexity} which is too expensive especially for large scale optimization problems.


A breakthrough on the quadratic termination property of the gradient method \eqref{eqitr} is due to Yuan \cite{yuan2006new,yuan2008step}, who obtained two-dimensional quadratic termination for the steepest descent (SD) method by deriving the so-called Yuan stepsize. More precisely, if the first and third steps use SD stepsizes while the second step uses the Yuan stepsize, then the minimizer of a two-dimensional strictly convex quadratic function is achieved in three iterations. By using a modified version of the Yuan stepsize and the SD stepsize in an alternate manner, Dai and Yuan \cite{dai2005analysis} developed the monotone Dai--Yuan method, which can even outperform the nonmonotone Barzilai--Borwein (BB) method \cite{Barzilai1988two}. As is known, the BB method, employing either of the following two genius stepsizes with certain quasi-Newton property:
\begin{equation}\label{bb1}
	\alpha_{k}^{BB1}=\arg \min _{\alpha>0}\left\|\alpha^{-1} s_{k-1}-y_{k-1}\right\|=\frac{s_{k-1}\tr s_{k-1}}{s_{k-1}\tr y_{k-1}}
\end{equation}
and
\begin{equation}\label{bb2}
	\alpha_{k}^{BB2}=\arg \min _{\alpha>0}\left\|s_{k-1}-\alpha y_{k-1}\right\|=\frac{s_{k-1}\tr y_{k-1}}{y_{k-1}\tr y_{k-1}},
\end{equation}
where $s_{k-1}=x_k-x_{k-1}$ and $y_{k-1}=g_k-g_{k-1}$, converges $R$-superlinearly for the two-dimensional strictly convex quadratic function \cite{Barzilai1988two} and performs much better than the SD method \cite{fletcher2005barzilai,raydan1997barzilai,yuan2008step}. It is noteworthy that Yuan's scheme of achieving two-dimensional quadratic termination is based on constructing a new coordinate system for $\mathbb{R}^2$ by two consecutive orthogonal gradients generated by the SD step, and transforming the solution from the standard rectangular coordinate system to the new one. Consequently, the Yuan stepsize and its variants  \cite{dai2005analysis,huang2022acceleration,huang2019asymptotic,li2023gradient} either depend on the use of the exact line search or  Hessian, which largely denies their extension to non-quadratic problems and many important applications.

A recent improvement in achieving two-dimensional quadratic termination for the gradient method \eqref{eqitr} is provided by Huang  et al. \cite{hdlbbq} who suggested to compute the stepsize from a necessary condition for the next gradient parallels to some eigenvector of the Hessian, see \eqref{stepmuprod2}. When applied to the BB method, the scheme in \cite{hdlbbq} yields
\begin{equation*}\label{snewa}
	\alpha_k^{BBQ}
	=\frac{2}{\frac{\phi_2}{\phi_3}+\sqrt{\left(\frac{\phi_2}{\phi_3}\right)^2-4\,\frac{\phi_1}{\phi_3}}},
\end{equation*}
where
\begin{equation*}\label{eqphinbb1}
	\frac{\phi_1}{\phi_3}
	=\frac{\alpha_{k-1}^{BB2}-\alpha_{k}^{BB2}}
	{\alpha_{k-1}^{BB2}\alpha_{k}^{BB2}(\alpha_{k-1}^{BB1}-\alpha_{k}^{BB1})}~~\mathrm{and}~~
	\frac{\phi_2}{\phi_3}
	=\frac{\alpha_{k-1}^{BB1}\alpha_{k-1}^{BB2}-\alpha_{k}^{BB1}\alpha_{k}^{BB2}}
	{\alpha_{k-1}^{BB2}\alpha_{k}^{BB2}(\alpha_{k-1}^{BB1}-\alpha_{k}^{BB1})}.
\end{equation*}
The stepsize $\alpha_k^{BBQ}$ not only ensures two-dimensional quadratic termination of the BB method with $\alpha_{k}^{BB1}$ or $\alpha_{k}^{BB2}$ but also has a great advantage of easy extension to general unconstrained optimization and a wide class of constrained optimization since its computation does not require exact line searches or the Hessian. Based on $\alpha_k^{BBQ}$ and the adaptive framework in \cite{zhou2006gradient}, an efficient gradient method for unconstrained optimization, namely BBQ, is developed in \cite{hdlbbq}, which performs significantly better than the BB method and other recent successful gradient methods. However, it is not clear whether  the gradient method \eqref{eqitr} with higher-dimensional quadratic termination leads further benefits.

In this paper, we introduce a simple and uniform mechanism for gradient methods to achieve three-dimensional quadratic termination. Based on the mechanism, a novel stepsize \eqref{un3ftalp} is derived such that the family of delayed gradient methods \eqref{fystepbb} equipped with the stepsize have the three-dimensional quadratic termination property, which immediately applies to the BB method. To our knowledge, this is the first result on three-dimensional quadratic termination of the BB method. A distinguished feature of the novel stepsize for the BB method is that, it does not require the use of any exact line search or the Hessian, and can be computed by stepsizes and gradient norms in recent three iterations.

Using long BB steps and some short steps associated with the novel stepsize in an adaptive manner, we develop an efficient gradient method, the method \eqref{snew}, for unconstrained quadratic optimization. Numerical experiments on quadratics show that the method \eqref{snew} performs much better than many successful gradient methods including BB \cite{Barzilai1988two}, Dai--Yuan (DY) \cite{yuan2008step},  SDC \cite{de2014efficient}, ABBmin2 \cite{frassoldati2008new}, and BBQ \cite{hdlbbq}. By making use of the Dai--Fletcher nonmonotone line search \cite{dai2005projected}, we further extend the method \eqref{snew} for unconstrained optimization which yields an efficient method, Algorithm \ref{alunc}. Numerical experiments on general unconstrained problems from the CUTEst collection \cite{gould2015cutest} demonstrate the advantage of Algorithm \ref{alunc} over BBQ \cite{hdlbbq}.

The paper is organized as follows. In Section \ref{sec2}, we introduce the mechanism for the gradient method to achieve three-dimensional quadratic termination which leads the derivation of a novel stepsize for the family \eqref{fystepbb}.  In Section \ref{secalgqp}, we specify the computation of the novel stepsize for the BB method and use it to develop an efficient gradient method, the method \eqref{snew}, for unconstrained quadratic optimization and an efficient method, Algorithm \ref{alunc}, for unconstrained optimization. Finally, some conclusion remarks are drawn in Section \ref{secclu}.

\section{A mechanism and a novel stepsize for the gradient method}\label{sec2}
In this section, we first introduce our mechanism of achieving three-dimensional quadratic termination for the gradient method \eqref{eqitr} by considering the unconstrained quadratic optimization
\begin{equation}\label{eqpro}
	\min_{x\in\mathbb{R}^n}~~f(x)=\frac{1}{2}x\tr Ax-b\tr x,
\end{equation}
where $b\in \mathbb{R}^n$ and $A\in \mathbb{R}^{n\times n}$ is symmetric positive definite. After presenting the mechanism, for a more general and uniform analysis, we utilize it to derive a novel stepsize for the family of delayed gradient methods with
\begin{equation}\label{fystepbb}
	\alpha_k=\frac{g_{k-1}\tr  \psi(A) g_{k-1}}{g_{k-1}\tr  \psi(A)A g_{k-1}},
\end{equation}
where $\psi$ is a real analytic function on $[\lambda_1,\lambda_n]$ with $\lambda_1$ and $\lambda_n$ being the smallest and largest eigenvalues of $A$, respectively, and can be expressed by Laurent series $\psi(z)=\sum_{k=-\infty}^\infty c_kz^k$ with $c_k\in\mathbb{R}$
such that $0<\sum_{k=-\infty}^\infty c_kz^k<+\infty$ for all $z\in[\lambda_1,\lambda_n]$. Clearly, the  two BB stepsizes $\alpha_k^{BB1}$ and $\alpha_k^{BB2}$ corresponding to $\psi(A)=I$ and $\psi(A)=A$, respectively.

\subsection{A mechanism for achieving three-dimensional quadratic termination}\label{secframework}
Consider the application of the gradient method \eqref{eqitr}  to the unconstrained quadratic optimization \eqref{eqpro} with $n=3$. Let $\mathrm{tr}(A)$ and $\mathrm{det}(A)$ be the trace and determinant of $A$, respectively, and $\lambda_1\leq\lambda_2\leq\lambda_3$ be the eigenvalues of $A$. Since 
\begin{equation*}\label{3dtrA}
	\mathrm{tr}(A)=\lambda_1+\lambda_2+\lambda_3
\end{equation*}
and
\begin{equation*}\label{3dtrAsq}
	\mathrm{tr}(A^2)=\lambda_1^2+\lambda_2^2+\lambda_3^2,
\end{equation*}
we have
\begin{equation*}
	\frac{\mathrm{tr}^2(A)-\mathrm{tr}(A^2)}{2}=\lambda_1\lambda_2+\lambda_1\lambda_3+\lambda_2\lambda_3.
\end{equation*}
Noticing that
\begin{equation*}\label{3ddetA}
	\mathrm{det}(A)=\lambda_1\lambda_2\lambda_3,
\end{equation*}
the roots of the cubic equation
\begin{equation}\label{qt3eq1}
	z^3-\mathrm{tr}(A)z^2+\frac{\mathrm{tr}^2(A)-\mathrm{tr}(A^2)}{2}z-\mathrm{det}(A)=0
\end{equation}
are exactly $\lambda_{1}$, $\lambda_{2}$ and $\lambda_{3}$. Thus a naive way to achieve three-dimensional quadratic termination is to use  reciprocals of the three roots of \eqref{qt3eq1} as three stepsizes. Of course this can be done for the three-dimensional case. The issue is whether we can use these stepsizes to design gradient methods for $n$-dimensional unconstrained quadratic optimization. Unfortunately, computing stepsizes from \eqref{qt3eq1} is not affordable for large $n$ since $\mathrm{det}(A)$ requires $\mathcal{O}(n^3)$ flops \cite{rezaifar2007new} and $\mathrm{tr}(A^2)$ takes $\mathcal{O}(n^2)$ flops. Moreover, the Hessian $A$ is not easy to obtain for a non-quadratic function which makes it difficult to extend these stepsizes for general unconstrained optimization. Hence, we have to achieve three-dimensional quadratic termination in a different way so that the stepsize can be computed in low cost and without involving the Hessian.

To proceed, we replace $A$ in \eqref{qt3eq1} by a matrix $H_k$ to get the following equation 
\begin{equation}\label{qt3eq4}
	z^3-\mathrm{tr}(H_k)z^2+\frac{\mathrm{tr}^2(H_k)-\mathrm{tr}(H_k^2)}{2}z-\mathrm{det}(H_k)=0,
\end{equation}
where
\begin{equation} \label{eqhesgft3}
	H_k
	=\begin{pmatrix}
		u_k\tr  A u_k & u_k\tr A v_k & u_k\tr A r_k\\
		v_k\tr A u_k   & v_k\tr A v_k & v_k\tr A r_k\\
		r_k\tr A u_k   & r_k\tr A v_k & r_k\tr A r_k\\
	\end{pmatrix}
\end{equation} 
with $u_k, v_k, r_k\in\mathbb{R}^n$. We take reciprocal of the largest root of the equation \eqref{qt3eq4} as our novel stepsize, namely $\alpha^{new}_k$. In the next section, we will see that, with proper $u_k, v_k, r_k$, the Hessian $A$ is not needed in computing $\alpha^{new}_k$ and all elements of $H_k$ can be calculated without additional matrix-vector products.

Before presenting our three-dimensional quadratic termination result, we need the following definition.
\begin{definition}
We say a stepsize $\alpha_{k}$ has the two-dimensional quadratic termination property if it satisfies
\begin{align}\label{stepmuprod2}
	&g_{\nu_1(k)} \tr \psi_1(A)(I-\alpha_kA)g_{k}\cdot g_{\nu_2(k)} \tr \psi_2(A)(I-\alpha_kA)g_{k}\nonumber\\
	&=g_{\nu_1(k)} \tr \psi_3(A)(I-\alpha_kA)g_{k}\cdot g_{\nu_2(k)} \tr \psi_4(A)(I-\alpha_kA)g_{k}
\end{align}
for some $\nu_1(k),\nu_2(k)\in\{1,\ldots,k\}$, where $\psi_1$, $\psi_2$, $\psi_3$, $\psi_4$ are real analytic functions on $[\lambda_1,\lambda_n]$ as $\psi$ in \eqref{fystepbb} and satisfy $\psi_1(z)\psi_2(z)=\psi_3(z)\psi_4(z)$ for $z\in\mathbb{R}$.
\end{definition}

Examples of stepsizes have the two-dimensional quadratic termination property including $\alpha_k^{BBQ}$, the Yuan stepsize \cite{yuan2006new} and its variants \cite{dai2005analysis,huang2022acceleration,huang2019asymptotic}, see \cite{hdlbbq} for more details. 

With the novel stepsize $\alpha^{new}_k$ and some stepsize has the two-dimensional quadratic termination property, we are able to achieve three-dimensional quadratic termination.
\begin{theorem}\label{3dftth}
	Consider the gradient method \eqref{eqitr} with the invariance property under orthogonal transformations for the problem \eqref{eqpro} for $n=3$. Then, if $u_k, v_k, r_k\in\mathbb{R}^3$ are orthonormal vectors and $\alpha_k=\alpha^{new}_k$, it holds that $g_{k+i}^{(3)}=0$ for all $i\geq1$. Further, if $\alpha_{k+j_0}$ has the two-dimensional quadratic termination property for some $j_0\geq3$, and $\alpha_{k+j}$ has the form 
	\begin{equation}\label{gnlfmalp1}
		\alpha_{k+j}=\left(\frac{g_{k+i_0}\tr  \psi(A) g_{k+i_0}}{g_{k+i_0}\tr  \psi(A)A^c g_{k+i_0}}\right)^{\frac{1}{c}},
	\end{equation}
for some $j\geq j_0+1$ with $i_0\in\{j_0+1,\ldots,j\}$, where $c$ is some real number and $\psi$ is a real analytic function on $[\lambda_1,\lambda_n]$ with the same property as the one in \eqref{fystepbb}, we must have $g_{k+i}=0$ for some $1\leq i\leq j+1$. 
\end{theorem}
\begin{proof}
	Since $u_k$, $v_k$ and $r_k$ are orthonormal vectors and $n=3$, $H_k$ and $A$ have the same eigenvalues due to the fact that $H_k=Q\tr AQ$ where $Q=(u_k,v_k,r_k)$ is an orthogonal matrix. Assume that the eigenvalues of $H_k$, say $\lambda_i(H_k)$, are arranged in increasing order, i.e., $\lambda_1(H_k)\leq\lambda_2(H_k)\leq\lambda_3(H_k)$. Notice that
	\begin{equation}\label{3dtrdeteqs2}
		\begin{aligned}
			\mathrm{tr}(H_k)&=\lambda_1(H_k)+\lambda_2(H_k)+\lambda_3(H_k),\\
			\mathrm{det}(H_k)&=\lambda_1(H_k)\lambda_2(H_k)\lambda_3(H_k),\\
			\mathrm{tr}(H_k^2)&=\lambda_1^2(H_k)+\lambda_2^2(H_k)+\lambda_3^2(H_k),
		\end{aligned}
	\end{equation}
	We must have $\alpha^{new}_k=1/\lambda_3(H_k)=1/\lambda_3$, which gives $g_{k+1}^{(3)}=(1-\alpha^{new}_k\lambda_3)g_{k}^{(3)}=0$ and $g_{k+i}^{(3)}=0$ for all $i\geq1$.

	If $g_{k+i}\neq0$, $i=1,\ldots,j_0$, and $\alpha_{k+j_0}$ has the two-dimensional quadratic termination property, by the proof of Theorem 2.1 in \cite{hdlbbq}, either $\alpha_{k+j_0}=1/\lambda_1$ or $\alpha_{k+j_0}=1/\lambda_2$. Without loss of generality, we assume $\alpha_{k+j_0}=1/\lambda_2$. Then $g_{k+i}$ has at most one nonzero component $g_{k+i}^{(1)}$ for all $i\geq j_0+1$. If $g_{k+i}\neq0$, $i=j_0+1,\ldots,j$,  and $\alpha_{k+j}$ has the form \eqref{gnlfmalp1}  for some $j\geq j_0+1$ with $i_0\in\{j_0+1,\ldots,j\}$, we have $\alpha_{k+j}=1/\lambda_1$ and thus $g_{k+j+1}=0$. This completes the proof. 
\end{proof}

From Theorem \ref{3dftth} we know that in order to achieve three-dimensional quadratic termination we need to take at least one gradient step with $\alpha^{new}_k$ and one gradient step with some stepsize $\alpha_k$ has the two-dimensional quadratic termination property. Moreover, some gradient steps are necessary before $\alpha^{new}_k$ and $\alpha_k$ so that the quantities involved in computing these two stepsizes are ready.	The form \eqref{gnlfmalp1} includes many popular stepsizes as special instances such as $\alpha^{BB1}_k$, $\alpha^{BB2}_k$, the SD stepsize $\alpha^{SD}_k$, and the stepsize $\alpha^{DAY}_k=\frac{\|s_{k-1}\|}{\|y_{k-1}\|}$ of the method in \cite{dai2015positive} (DAY for short). Thus, three-dimensional quadratic termination will be achieved for the aforementioned methods if we take any of them for some steps before/after $\alpha^{new}_k$ and $\alpha_k$.


To numerical verify the results of three-dimensional quadratic termination, for the DAY, BB1 and BB2 methods, we applied the above described steps to the unconstrained quadratic optimization problem \eqref{eqpro} with
\begin{equation}\label{twoquad}
	A=\textrm{diag}\{1,\,\kappa/2,\,\kappa\} \quad \mbox{and} \quad b=0.
\end{equation}
In particular, we used $\alpha_1=\alpha^{SD}_1$, $\alpha_3=\alpha^{new}_3$ and  $\alpha_6=\alpha^{BBQ}_6$ for all three methods where $\alpha^{new}_3$ is computed with $u_3,v_3,r_3$ obtained by applying the Gram--Schmidt orthogonalization to $g_1,g_2,g_3$. When $k\neq 1,3,6$, we used $\alpha_k=\alpha^{DAY}_k$,  $\alpha_k=\alpha^{BB1}_k$ and $\alpha_k=\alpha^{BB2}_k$ for the DAY, BB1 and BB2 methods, respectively. 
By the proof of Theorem \ref{3dftth}, these three methods, referred to as DAY3D, BB1-3D, and BB2-3D for short, respectively, will give the solution at $x_9$. 
Table \ref{tb3ft} presents averaged values of gradient norm and objective value over ten random starting points. We can see that for different $\kappa$, the gradient norms and objective values obtained by DAY3D, BB1-3D, and BB2-3D are numerically very close to zero while those values obtained by the unmodified BB1 method are far away from zero.

\begin{table}[ht!b]
	\setlength{\tabcolsep}{0.98ex}
	{\footnotesize
		\caption{Averaged results on problem \eqref{twoquad} with different condition numbers.}\label{tb3ft}
		\begin{center}
			\begin{tabular}{|c|c|c|c|c|c|c|c|c|}
				\hline
				\multicolumn{1}{|c|}{\multirow{2}{*}{$\kappa$}}  &\multicolumn{2}{c|}{\multirow{1}{*}{BB1}}   &\multicolumn{2}{c|}{\multirow{1}{*}{DAY3D}}  &\multicolumn{2}{c|}{\multirow{1}{*}{BB1-3D}}  &\multicolumn{2}{c|}{\multirow{1}{*}{BB2-3D}}  \\
				\cline{2-9}
				&\multicolumn{1}{c|}{$\|g_9\|$} &\multicolumn{1}{c|}{$f(x_9)$}   &\multicolumn{1}{c|}{$\|g_9\|$} &\multicolumn{1}{c|}{$f(x_9)$}  &\multicolumn{1}{c|}{$\|g_9\|$} &\multicolumn{1}{c|}{$f(x_9)$}  &\multicolumn{1}{c|}{$\|g_9\|$} &\multicolumn{1}{c|}{$f(x_9)$}  \\
				\hline
100    &2.80e+00  &2.24e-01  &1.59e-12  &3.15e-26  &3.12e-12  &1.25e-25  &1.53e-12  &2.14e-26\\
\hline
1000   &7.12e-01  &1.65e-01  &2.31e-10  &1.27e-22  &3.64e-10  &3.11e-22  &2.07e-10  &7.08e-23\\
\hline
10000  &1.94e+01  &3.66e-01  &4.06e-08  &5.62e-19  &1.12e-07  &4.59e-18  &1.41e-07  &7.72e-18\\
				\hline
			\end{tabular}
		\end{center}
	}
\end{table}

\subsection{A novel stepsize equipped for the family \eqref{fystepbb}}\label{secBBSD}

We first give the formula of the novel stepsize $\alpha^{new}_k$ in the following lemma.
\begin{lemma}\label{3alpdefth}
	Assume that $u_k, v_k, r_k$ are not all null vectors. Then $\alpha^{new}_k$ has the form
	\begin{equation}\label{un3ftalp}
		\alpha^{new}_k=\frac{1}{\frac{\mathrm{tr}(H_k)}{3}+2\cos \left(\frac{\theta_k}{3}\right)\sqrt{\frac{|p_k|}{3}}},
	\end{equation}
	where
	\begin{equation*}
		\theta_k=\arccos\left(-\frac{q_k}{2}\left(\frac{3}{|p_k|}\right)^{3/2}\right),
		\quad
		p_k=\frac{\mathrm{tr}^2(H_k)-3\mathrm{tr}(H_k^2)}{6},
	\end{equation*}
	and
	\begin{equation*}
		q_k=\frac{5\mathrm{tr}^3(H_k)-9\mathrm{tr}(H_k)\mathrm{tr}(H_k^2)}{54}-\mathrm{det}(H_k).
	\end{equation*}	
\end{lemma}
\begin{proof}
	The formula of $\alpha^{new}_k$ can be obtained by standard deduction of roots for cubic equations, see \cite{zucker200892} for example. 
\end{proof}


Interestingly, the stepsize $\alpha^{new}_{k}$ is bounded by reciprocals of the trace and diagonal elements of $H_k$ even for the general $n$-dimensional case.
\begin{theorem}\label{bdalpnew}
	If $u_k, v_k$ and $r_k$ are orthonormal vectors,  then	
	\begin{align}\label{un3ftalpbdlm}
		\frac{1}{\mathrm{tr}(H_k)}\leq
		&\alpha^{new}_k\leq\min\left\{\frac{1}{H_k^{(11)}},\frac{1}{H_k^{(22)}},\frac{1}{H_k^{(33)}}\right\},
	\end{align}
	where $H_k^{(ij)}$ is the $(i,j)$ component of $H_k$.
\end{theorem}
\begin{proof}
	Let $\hat{H}_k$ be a principal submatrix of $H_k$ obtained by deleting one row and the corresponding column. Denote $\lambda_1(\hat{H}_k)\leq\lambda_2(\hat{H}_k)\leq\lambda_3(\hat{H}_k)$ as the eigenvalues of $\hat{H}_k$. Since $u_k, v_k$ and $r_k$ are orthonormal vectors, by Theorem 4.3.15 of  \cite{horn2012matrix}, we have
	\begin{equation}\label{eigsrel1}
		\lambda_1(H_k)\leq\lambda_1(\hat{H}_k)\leq\lambda_2(H_k)\leq\lambda_2(\hat{H}_k)\leq\lambda_3(H_k).
	\end{equation}
	Let
	\begin{equation*}\label{qt32dalp}
		\alpha_k^{\hat{H}}
		=\frac{2}{\mathrm{tr}(\hat{H}_k)+\sqrt{\mathrm{tr}^2(\hat{H}_k)-4\mathrm{det}(\hat{H}_k)}}.
	\end{equation*} 
	Similar to \eqref{3dtrdeteqs2}, we get $\alpha_k^{\hat{H}}=1/\lambda_2(\hat{H}_k)$,	which together with \eqref{eigsrel1} gives $\alpha^{new}_k=1/\lambda_3(H_k)\leq\alpha_k^{\hat{H}}$. The right inequality of \eqref{un3ftalpbdlm} follows by the fact that $\alpha_k^{\hat{H}}\leq\frac{1}{\max\{\hat{H}_k^{(11)},\hat{H}_k^{(22)}\}}$ and the consideration of different submatrices.

	Since $0\leq\theta_k\leq\pi$, it holds that $1/2\leq\cos \left(\theta_k/3\right)\leq1$, which together with \eqref{un3ftalp} indicates that
	\begin{align*}\label{un3ftalpbd1}
		\alpha^{new}_k\geq\frac{1}{\frac{\mathrm{tr}(H_k)}{3}+2\sqrt{\frac{|p_k|}{3}}}
		=\frac{3}{\mathrm{tr}(H_k)+2\sqrt{\frac{3\mathrm{tr}(H_k^2)-\mathrm{tr}^2(H_k)}{2}}}.
	\end{align*}	
	To prove the left inequality of \eqref{un3ftalpbdlm}, we are suffice to show
	\begin{equation*}
		\mathrm{tr}(H_k)\geq\sqrt{\frac{3\mathrm{tr}(H_k^2)-\mathrm{tr}^2(H_k)}{2}},
	\end{equation*}
	which is clearly true since $\mathrm{tr}^2(H_k)\geq\mathrm{tr}(H_k^2)$.  This completes the proof.
\end{proof}

A pre-requisite for computing $\alpha^{new}_k$ for the family \eqref{fystepbb} is choosing orthonormal vectors $u_k$, $v_k$ and $r_k$ to determine the matrix $H_k$. There are various ways to obtain orthonormal vectors such as Gram--Schmidt process, Householder transformation, Givens rotation, see \cite{golub2013matrix} for example. Here we apply the Gram--Schmidt process to historical gradients which leads to easy computation of the matrix $H_k$. In particular, we set
\begin{equation}\label{uvk3ftbb}
	u_k=\frac{\psi^{c}(A)g_{k-3}}{\|\psi^{c}(A)g_{k-3}\|},
	\quad
	v_k=\frac{\bar{v}_k}{\|\bar{v}_k\|},\quad
	r_k=\frac{\bar{r}_k}{\|\bar{r}_k\|},
\end{equation}
with
\begin{equation*}
	\bar{v}_k=\psi^{1-c}(A)g_{k-2}-\frac{g_{k-2}\tr \psi(A)g_{k-3}}{\|\psi^c(A)g_{k-3}\|^2}\psi^c(A)g_{k-3}
\end{equation*}
and
\begin{equation*}
	\bar{r}_k=\psi^c(A)g_{k-1}-\frac{g_{k-1}\tr \psi^{2c}(A)g_{k-3}}{\|\psi^c(A)g_{k-3}\|^2}\psi^c(A)g_{k-3}-\frac{g_{k-1}\tr \psi^c(A)\bar{v}_k}{\|\bar{v}_k\|^2}\bar{v}_k,
\end{equation*}
where $c\in\mathbb{R}$. Under the assumption that $g_{k-1}$, $g_{k-2}$ and $g_{k-3}$ are linearly independent, $u_k$, $v_k$ and $r_k$ are orthonormal.

%
%
%
%

Now we specify the matrix $H_k$ for the family \eqref{fystepbb} with $u_k$, $v_k$ and $r_k$ given by \eqref{uvk3ftbb}. To distinguish different stepsizes, we refer $\alpha_{k}$  defined by  \eqref{fystepbb} as $\alpha_{k}^\psi$.
\begin{theorem}\label{thhesbb12}
	Consider the family \eqref{fystepbb} for the problem \eqref{eqpro}. Suppose  $g_{k-1}$, $g_{k-2}$ and $g_{k-3}$ are linearly independent, and $u_k$, $v_k$ and $r_k$ are given by \eqref{uvk3ftbb}. For $c=1/2$, the matrix $H_{k}$ has the form
	\begin{equation}\label{eqhesbb12ft3}
		H_{k}=\begin{pNiceMatrix}
			\Block{2-2}{\hat{H}_k} & & 0\\
			\hspace*{1cm} & &-\frac{\sqrt{g_{k-1}\tr \psi^{1/2}(A)\bar{r}_k}}{\alpha_{k-2}\|\bar{v}_k\|} \\
			0   &
			-\frac{\sqrt{g_{k-1}\tr \psi^{1/2}(A)\bar{r}_k}}{\alpha_{k-2}\|\bar{v}_k\|}
			& 	\frac{g_{k-1}\tr \psi^{1/2}(A)A\bar{r}_k}{g_{k-1}\tr \psi^{1/2}(A)\bar{r}_k}+\frac{\gamma_k}{\alpha_{k-2}}
		\end{pNiceMatrix}
	\end{equation}
	with
	\begin{equation}\label{eqhathesmonbb}
		\hat{H}_k= \begin{pmatrix}
			\frac{1}{\alpha_{k-2}^\psi} & -\frac{\sqrt{1-\sigma_k}\|\psi^{1/2}(A)g_{k-2}\|}{\alpha_{k-3}\|\psi^{1/2}(A)g_{k-3}\|}\\
			-\frac{\sqrt{1-\sigma_k}\|\psi^{1/2}(A)g_{k-2}\|}{\alpha_{k-3}\|\psi^{1/2}(A)g_{k-3}\|}  & \frac{1}{1-\sigma_k}\left(\frac{1}{\alpha_{k-1}^\psi}-
			2\sigma_k\delta_k+\frac{\sigma_k}{\alpha_{k-2}^\psi}\right)
		\end{pmatrix},
	\end{equation}	
	where 
	\begin{equation*}
		\|\bar{v}_k\|=\|\psi^{1/2}(A)g_{k-2}\|\sqrt{1-\sigma_k},
		\quad
		\gamma_k=\frac{g_{k-1}\tr \psi^{1/2}(A)\bar{v}_k}{\|\bar{v}_k\|^2},
	\end{equation*}
	\begin{equation*}
		\sigma_k=\frac{(g_{k-2}\tr \psi(A)g_{k-3})^2}{\|\psi^{1/2}(A)g_{k-3}\|^2\|\psi^{1/2}(A)g_{k-2}\|^2},\quad\delta_k=\frac{g_{k-2}\tr \psi(A) Ag_{k-3}}{g_{k-2}\tr \psi(A)g_{k-3}}.
	\end{equation*}
\end{theorem}
\begin{proof}
	By the definitions of $u_k$, $\bar{v}_k$ and $\bar{r}_k$, we have
	\begin{equation}\label{orthbb12}
		g_{k-3}\tr \psi^{1/2}(A)\bar{v}_k=0,\quad 
		g_{k-3}\tr \psi^{1/2}(A)\bar{r}_k=0,\quad 
		\bar{v}_k\tr \bar{r}_k=0,
	\end{equation}
	\begin{align}\label{uAubb12}
		u_k\tr  Au_k&=\frac{g_{k-3}\tr \psi(A)Ag_{k-3}}{\|\psi^{1/2}(A)g_{k-3}\|^2}=\frac{1}{\alpha_{k-2}^\psi},
	\end{align}	
	and
	\begin{align}\label{bvavbb12}
		\bar{v}_k\tr A\bar{v}_k
		&=g_{k-2}\tr  \psi(A)Ag_{k-2}-2\frac{g_{k-2}\tr \psi(A)g_{k-3}}{\|\psi^{1/2}(A)g_{k-3}\|^2}g_{k-2}\tr  \psi(A)Ag_{k-3}\nonumber\\
		&\quad\quad\quad+\frac{(g_{k-2}\tr \psi(A)g_{k-3})^2}{\|\psi^{1/2}(A)g_{k-3}\|^4}g_{k-3}\tr \psi(A) Ag_{k-3}\nonumber\\
		&=\|\psi^{1/2}(A)g_{k-2}\|^2\left(\frac{1}{\alpha_{k-1}^\psi}-2\sigma_k\delta_k+\frac{\sigma_k}{\alpha_{k-2}^\psi}\right).
	\end{align}	
	From \eqref{orthbb12}, we obtain
	\begin{equation}\label{orthgk2rbb12}
		g_{k-2}\tr \psi^{1/2}(A)\bar{r}_k=0,
	\end{equation}
	\begin{align}\label{rrbb12}
		\bar{r}_k\tr \bar{r}_k
		&=g_{k-1}\tr \psi^{1/2}(A)\bar{r}_k
	\end{align}
	and
	\begin{align}\label{vvbb12}
		\|\bar{v}_k\|^2&=
		g_{k-2}\tr \psi^{1/2}(A)\bar{v}_k
		=\|\psi^{1/2}(A)g_{k-2}\|^2\left(1-\sigma_k\right).
	\end{align}
	Combining \eqref{bvavbb12} and  \eqref{vvbb12}, we get 
	\begin{align}\label{vAvbb12}
		v_k\tr  Av_k&=\frac{\bar{v}_k\tr A\bar{v}_k}{\|\bar{v}_k\|^2}=\frac{1}{1-\sigma_k}
		\left(\frac{1}{\alpha_{k-1}^\psi}-2\sigma_k\delta_k+\frac{\sigma_k}{\alpha_{k-2}^\psi}\right).
	\end{align}
	Noting that
	\begin{equation}\label{Agkbb12}
		Ag_{k-j}=\frac{1}{\alpha_{k-j}}(g_{k-j}-g_{k-j+1}),~~j=0,1,2,\ldots,
	\end{equation}
	which together with \eqref{orthbb12} and \eqref{vvbb12} yields
	\begin{align}\label{uAvbb12}
		u_k\tr  Av_k&=\frac{(g_{k-3}-g_{k-2})\tr \psi^{1/2}(A) \bar{v}_k}{\alpha_{k-3}\|\psi^{1/2}(A)g_{k-3}\|\|\bar{v}_k\|}
		=-\frac{g_{k-2}\tr \psi^{1/2}(A)\bar{v}_k}{\alpha_{k-3}\|\psi^{1/2}(A)g_{k-3}\|\|\bar{v}_k\|}
		\nonumber\\
		&=-\frac{\|\bar{v}_k\|}{\alpha_{k-3}\|\psi^{1/2}(A)g_{k-3}\|}
		=-\frac{\sqrt{1-\sigma_k}\|\psi^{1/2}(A)g_{k-2}\|}{\alpha_{k-3}\|\psi^{1/2}(A)g_{k-3}\|}.
	\end{align}	
	By \eqref{orthbb12},  \eqref{orthgk2rbb12} and \eqref{Agkbb12}, we obtain
	\begin{align}\label{gk3Arbb12}
		g_{k-3}\tr \psi^{1/2}(A)A\bar{r}_k&=\frac{1}{\alpha_{k-3}}(g_{k-3}-g_{k-2})\tr \psi^{1/2}(A)\bar{r}_k=0,
	\end{align}
	which indicates
	\begin{align}\label{uArbb12}
		u_{k}\tr Ar_k&=0
	\end{align}
	and
	\begin{align}\label{bvArbb12}
		\bar{v}_k\tr A\bar{r}_k
		&=g_{k-2}\tr \psi^{1/2}(A)A\bar{r}_k-\frac{g_{k-2}\tr \psi(A)g_{k-3}}{\|\psi^{1/2}(A)g_{k-3}\|^2}g_{k-3}\tr \psi^{1/2}(A)A\bar{r}_k
		\nonumber\\
		&=\frac{1}{\alpha_{k-2}}(g_{k-2}-g_{k-1})\tr \psi^{1/2}(A)\bar{r}_k=-\frac{g_{k-1}\tr \psi^{1/2}(A)\bar{r}_k}{\alpha_{k-2}}.
	\end{align}
	Then, by \eqref{rrbb12} and \eqref{bvArbb12}, we have
	\begin{align}\label{vArbb12}
		v_{k}\tr Ar_k&=\frac{\bar{v}_k\tr A\bar{r}_k}{\|\bar{v}_k\|\|\bar{r}_k\|}=-\frac{\sqrt{g_{k-1}\tr \psi^{1/2}(A)\bar{r}_k}}{\alpha_{k-2}\|\bar{v}_k\|}.
	\end{align}
	From \eqref{gk3Arbb12} and \eqref{bvArbb12}, we obtain
	\begin{align*}
		\bar{r}_k\tr A\bar{r}_k
		&=g_{k-1}\tr \psi^{1/2}(A)A\bar{r}_k+\frac{g_{k-1}\tr \psi^{1/2}(A)\bar{v}_k}{\|\bar{v}_k\|^2}\frac{1}{\alpha_{k-2}}g_{k-1}\tr \psi^{1/2}(A)\bar{r}_k\\
		&=g_{k-1}\tr \psi^{1/2}(A)A\bar{r}_k+\frac{g_{k-1}\tr \psi^{1/2}(A)\bar{r}_k}{\alpha_{k-2}}\gamma_k,
	\end{align*}			
	which combines \eqref{uAubb12}, \eqref{vvbb12}, \eqref{vAvbb12}, \eqref{uAvbb12}, \eqref{uArbb12}, and \eqref{vArbb12}   gives \eqref{eqhesbb12ft3}. This completes the proof.
\end{proof}

\begin{remark}
	We note that the limited memory steepest descent (LMSD) method proposed by Fletcher \cite{fletcher2012limited} also employs recent gradients to determine the stepsize, which is proved to be $R$-linear convergent for strongly convex quadratics \cite{Curtis2018}. However, it needs to keep some gradients in memory and apply thin QR factorization or extended Cholesky factorization to the matrix associated with these gradients to compute stepsizes. In the next section, we will see that, for the BB method with $\alpha_{k}^{BB1}$, i.e. $\psi(A)=I$, the stepsize $\alpha^{new}_{k}$ can be computed by using stepsizes and gradient norms in recent three iterations. Hence we only need to keep these scalars in memory without conducting factorization on a matrix, nor requiring additional matrix-vector and inner products.
\end{remark}

\section{Application to the BB method}\label{secalgqp}
In this section, we apply the proposed mechanism to the BB method and develop an efficient gradient method for solving unconstrained optimization problems.

We will focus on the BB method with $\alpha_{k}^{BB1}$ since the analysis also applies to that with $\alpha_{k}^{BB2}$ and practical experience often prefers $\alpha_{k}^{BB1}$ (see \cite{dai2005asymptotic}). We proceed by specifying the matrix $H_{k}$ for the BB method. Interestingly, it only depends on stepsizes and gradient norms in recent three iterations, and does not require exact line searches or the Hessian.

\begin{corollary}\label{corhesbb1}
	Under the conditions of Theorem \ref{thhesbb12}, for the BB method, the matrix $H_{k}$ has the form
	\begin{equation*}\label{hmfBB}
		H_{k}=\begin{pNiceMatrix}
			\Block{2-2}{\hat{H}_k} & & 0\\
			\hspace*{1cm} & & -\frac{\sqrt{g_{k-1}\tr \bar{r}_k}}{\alpha_{k-2}\|g_{k-2}\|\sqrt{1-\sigma_k}}
			\\
			0   & -\frac{\sqrt{g_{k-1}\tr \bar{r}_k}}{\alpha_{k-2}\|g_{k-2}\|\sqrt{1-\sigma_k}}
			& 
			\frac{g_{k-1}\tr A\bar{r}_k}{g_{k-1}\tr \bar{r}_k}
			+\frac{\gamma_k}{\alpha_{k-2}}
		\end{pNiceMatrix}
	\end{equation*}
	with
	\begin{equation*}
		\hat{H}_k= \begin{pmatrix}
			\frac{1}{\alpha_{k-2}^{BB1}} & -\frac{\sqrt{1-\sigma_k}\|g_{k-2}\|}{\alpha_{k-3}\|g_{k-3}\|}\\
			-\frac{\sqrt{1-\sigma_k}\|g_{k-2}\|}{\alpha_{k-3}\|g_{k-3}\|}  & \frac{1}{1-\sigma_k}\left(\frac{1}{\alpha_{k-1}^{BB1}}-
			2\sigma_k\delta_k+\frac{\sigma_k}{\alpha_{k-2}^{BB1}}\right)
		\end{pmatrix},
	\end{equation*}	
	where 
	\begin{equation*}
		\gamma_k=1-\frac{\alpha_{k-2}}{1-\sigma_k}\left(\frac{1}{\alpha_{k-1}^{BB1}}-\sigma_k\delta_k\right),
		\quad
		\sigma_k
		=\left(1-\frac{\alpha_{k-3}}{\alpha_{k-2}^{BB1}}\right)\zeta_k,
	\end{equation*}
	and
	\begin{equation*}
		\delta_k=\frac{1}{\alpha_{k-3}}\left(1-\frac{1}{\zeta_k}\right),
		\quad
		\zeta_k=\left(1-\frac{\alpha_{k-3}}{\alpha_{k-2}^{BB1}}\right)\frac{\|g_{k-3}\|^2}{\|g_{k-2}\|^2}.
	\end{equation*}
	In addition,
	\begin{align*}
		g_{k-1}\tr \bar{r}_k
		&=\|g_{k-1}\|^2-\left(\sigma_k(1-\alpha_{k-2}\delta_k)^2+\gamma_k^2(1-\sigma_k)\right)\|g_{k-2}\|^2
	\end{align*}	
	and
	\begin{align*}
		g_{k-1}\tr A\bar{r}_k
		&=\left(\frac{1}{\alpha_{k}^{BB1}}+\frac{\gamma_k}{\alpha_{k-2}}\right)\|g_{k-1}\|^2+\varsigma_k\|g_{k-2}\|^2,
	\end{align*}
	where
	\begin{align*}
		\varsigma_k&=\left(\frac{\gamma_k-\left(1-\alpha_{k-2}\delta_k\right)}{\alpha_{k-2}^{BB1}}-\frac{\gamma_k}{\alpha_{k-2}}\right)\left(1-\frac{\alpha_{k-2}}{\alpha_{k-1}^{BB1}}\right)\\
		&\quad
		-\frac{\gamma_k-\left(1-\alpha_{k-2}\delta_k\right)}{\alpha_{k-3}}\gamma_k(1-\sigma_k).
	\end{align*}	
\end{corollary}	
\begin{proof}
	From Theorem \ref{thhesbb12}, \eqref{Agkbb12} and
	\begin{equation}\label{gktgk1}
		g_{k-j+1}\tr g_{k-j}=\left(1-\frac{\alpha_{k-j}}{\alpha_{k-j+1}^{BB1}}\right)\|g_{k-j}\|^2,~~j=0,1,2,\ldots,
	\end{equation}
	we get
	\begin{equation*}
		\sigma_k
		=\frac{(g_{k-2}\tr g_{k-3})^2}{\|g_{k-2}\|^2\|g_{k-3}\|^2}
		=\left(1-\frac{\alpha_{k-3}}{\alpha_{k-2}^{BB1}}\right)^2\frac{\|g_{k-3}\|^2}{\|g_{k-2}\|^2}
		=\left(1-\frac{\alpha_{k-3}}{\alpha_{k-2}^{BB1}}\right)\zeta_k,
	\end{equation*}
	and
	\begin{equation*}
		\delta_k=\frac{g_{k-2}\tr Ag_{k-3}}{g_{k-2}\tr g_{k-3}}=\frac{1}{\alpha_{k-3}}
		\left(1-\frac{\|g_{k-2}\|^2}{\|g_{k-3}\|^2\left(1-\frac{\alpha_{k-3}}{\alpha_{k-2}^{BB1}}\right)}\right)
		=\frac{1}{\alpha_{k-3}}\left(1-\frac{1}{\zeta_k}\right).
	\end{equation*}
	By $g_{k-1}\tr g_{k-3}=g_{k-2}\tr g_{k-3}(1-\alpha_{k-2}\delta_k)$, \eqref{gktgk1} and the definition of $\sigma_k$, we have
	\begin{align*}
		g_{k-1}\tr \bar{v}_k&=g_{k-1}\tr g_{k-2}-\frac{g_{k-2}\tr g_{k-3}}{\|g_{k-3}\|^2}g_{k-1}\tr g_{k-3}\nonumber\\
		&=g_{k-1}\tr g_{k-2}-\sigma_k(1-\alpha_{k-2}\delta_k)\|g_{k-2}\|^2\nonumber\\
		&=\left(1-\sigma_k-\alpha_{k-2}\left(\frac{1}{\alpha_{k-1}^{BB1}}-\sigma_k\delta_k\right)\right)\|g_{k-2}\|^2,
	\end{align*}
	which together with the fact $\|\bar{v}_k\|^2=\|g_{k-2}\|^2\left(1-\sigma_k\right)$ yields the formula of $\gamma_k$. 
	
	By  the definitions of $\sigma_k$ $\bar{r}_k$, and $\bar{v}_k$, and again $g_{k-1}\tr g_{k-3}=g_{k-2}\tr g_{k-3}(1-\alpha_{k-2}\delta_k)$ and $\|\bar{v}_k\|^2=\|g_{k-2}\|^2\left(1-\sigma_k\right)$, we obtain
	\begin{align*}
		g_{k-1}\tr \bar{r}_k
		&=\|g_{k-1}\|^2-\frac{(g_{k-1}\tr g_{k-3})^2}{\|g_{k-3}\|^2}-\gamma_k^2\|\bar{v}_k\|^2\\
		&=\|g_{k-1}\|^2-\left(\sigma_k(1-\alpha_{k-2}\delta_k)^2+\gamma_k^2(1-\sigma_k)\right)\|g_{k-2}\|^2
	\end{align*}
	and
	\begin{align*}
		\bar{r}_k
		&=g_{k-1}-\frac{g_{k-2}\tr g_{k-3}\left(1-\alpha_{k-2}\delta_k\right)}{\|g_{k-3}\|^2}g_{k-3}-\gamma_k\bar{v}_k\\
		&=g_{k-1}-\gamma_kg_{k-2}+
		\frac{g_{k-2}\tr g_{k-3}}{\|g_{k-3}\|^2}\left(\gamma_k
		-1+\alpha_{k-2}\delta_k
		\right)g_{k-3}.
	\end{align*}
	Then, by \eqref{Agkbb12} and \eqref{gktgk1}, we get
	\begin{align*}
		g_{k-1}\tr A\bar{r}_k
		&=\frac{\|g_{k-1}\|^2}{\alpha_{k}^{BB1}}-\frac{\gamma_k}{\alpha_{k-2}}(g_{k-1}\tr g_{k-2}-\|g_{k-1}\|^2)\\
		&\quad\quad +
		\frac{g_{k-2}\tr g_{k-3}}{\|g_{k-3}\|^2}\frac{\gamma_k-\left(1-\alpha_{k-2}\delta_k\right)}{\alpha_{k-3}}
		\left(g_{k-1}\tr g_{k-3}-g_{k-1}\tr g_{k-2}\right)\\
		&=\left(\frac{1}{\alpha_{k}^{BB1}}+\frac{\gamma_k}{\alpha_{k-2}}\right)\|g_{k-1}\|^2-\frac{\gamma_k}{\alpha_{k-2}}g_{k-1}\tr g_{k-2}+\\
		&\quad\quad  \frac{\gamma_k-\left(1-\alpha_{k-2}\delta_k\right)}{\alpha_{k-3}}
		\left(g_{k-1}\tr g_{k-2}-g_{k-1}\tr \bar{v}_k-\frac{g_{k-2}\tr g_{k-3}}{\|g_{k-3}\|^2}g_{k-1}\tr g_{k-2}\right)
		\\
		&=\left(\frac{1}{\alpha_{k}^{BB1}}+\frac{\gamma_k}{\alpha_{k-2}}\right)\|g_{k-1}\|^2-\frac{\gamma_k}{\alpha_{k-2}}g_{k-1}\tr g_{k-2}+\\
		&\quad\quad  \frac{\gamma_k-\left(1-\alpha_{k-2}\delta_k\right)}{\alpha_{k-3}}
		\left(\frac{\alpha_{k-3}}{\alpha_{k-2}^{BB1}}g_{k-1}\tr g_{k-2}-g_{k-1}\tr \bar{v}_k\right)
		\\
		&=\left(\frac{1}{\alpha_{k}^{BB1}}+\frac{\gamma_k}{\alpha_{k-2}}\right)\|g_{k-1}\|^2+\left(\frac{\gamma_k-\left(1-\alpha_{k-2}\delta_k\right)}{\alpha_{k-2}^{BB1}}-\frac{\gamma_k}{\alpha_{k-2}}\right)g_{k-1}\tr g_{k-2}\\
		&\quad\quad -\frac{\gamma_k-\left(1-\alpha_{k-2}\delta_k\right)}{\alpha_{k-3}}
		\gamma_k(1-\sigma_k)\|g_{k-2}\|^2\\
		&=\left(\frac{1}{\alpha_{k}^{BB1}}+\frac{\gamma_k}{\alpha_{k-2}}\right)\|g_{k-1}\|^2+\varsigma_k\|g_{k-2}\|^2.
	\end{align*}
	This completes the proof.
\end{proof}

\subsection{Quadratic optimization}
In this subsection, we concentrate on the quadratic optimization \eqref{eqpro} the study of which is  usually a pre-requisite for developing and analyzing gradient methods for general nonlinear optimization.

It is commonly accepted that incorporating the long BB stepsize $\alpha_{k}^{BB1}$ (since $\alpha_{k}^{BB1}\geq\alpha_{k}^{BB2}$) with some short stepsizes under the adaptive scheme \cite{Bonettini2008,hdlbbq,zhou2006gradient} often performs better than the original BB method, see \cite{Serena2020,dai2006cyclic,frassoldati2008new,hdlbbq,huang2022acceleration,zhou2006gradient} for example. 
From Theorem \ref{bdalpnew} and  Corollary \ref{corhesbb1}, we get
\begin{align*}\label{un3ftalpupbdbb1}
	&\alpha^{new}_k\leq\min\left\{\frac{1}{H_k^{(11)}},\frac{1}{H_k^{(22)}},\frac{1}{H_k^{(33)}}\right\}\leq\alpha_{k-2}^{BB1}.
\end{align*}
In this sense, the stepsize $\alpha^{new}_k$ is short. So, we would like to develop a gradient method that adaptively taking the long BB stepsize $\alpha_{k}^{BB1}$ and some short stepsize associated with $\alpha^{new}_{k}$.

However, we do not know whether $\alpha_k^{new}$ is larger than the recent short BB stepsizes. To ensure a short stepsize, we shall utilize $\min\{\alpha_{k-1}^{BB2},\alpha_{k}^{BB2},\alpha_k^{new}\}$ instead of $\alpha_k^{new}$. Due to rounding error and loss of  linear independence of historical gradients, the term $g_{k-1}\tr \bar{r}_k$ may be nonpositive. In this case, $\alpha_k^{new}$ is not well-defined and we  replace it by $\alpha_k^{BBQ}$. To sum up,
our method employs the stepsizes $\{\alpha_k: k\ge 5\}$ given by 
\begin{equation}\label{snew}
	\alpha_{k}= 
	\begin{cases}
		\min\{\alpha_{k-1}^{BB2},\alpha_{k}^{BB2},\alpha_k^{new}\}, & \hbox{if $\alpha_k^{BB2}/\alpha_k^{BB1}<\tau_k$ and $g_{k-1}\tr \bar{r}_k>0$;} \\
		\min\{\alpha_{k-1}^{BB2},\alpha_{k}^{BB2},\alpha_k^{BBQ}\}, & \hbox{if $\alpha_k^{BB2}/\alpha_k^{BB1}<\tau_k$ and $g_{k-1}\tr \bar{r}_k\leq0$;} \\
		\alpha_k^{BB1}, & \hbox{otherwise,}
	\end{cases}
\end{equation}
where $\alpha^{new}_k$ is computed with $H_k$ given in Corollary \ref{corhesbb1} and $\tau_k>0$ is updated by
\begin{equation} \label{tauk}
	\tau_{k+1}=\left\{
	\begin{array}{ll}
		\tau_k/\gamma, & \hbox{ if $\alpha_k^{BB2}/\alpha_k^{BB1}<\tau_k$;} \\
		\tau_k\gamma, & \hbox{ otherwise,}
	\end{array}
	\right.
\end{equation} 
for some $\gamma\geq1$, see \cite{Bonettini2008,hdlbbq} for example. Moreover, we take $\alpha_1=\alpha_1^{SD}$ and $\alpha_i=\alpha_i^{BB1}$, $i=2,3,4$, for quadratic optimization \eqref{eqpro}.

To establish $R$-linear convergence of the method \eqref{snew}, we show that $\alpha_k$ in \eqref{snew} satisfies the Property B in \cite{li2023note}. In particular, we say that a stepsize $\alpha_k$ has Property B if there exist an integer $m$ and a positive constant $M_{1} \geq \lambda_{1}$ such that\\
(i) $\lambda_{1} \leq \alpha_{k}^{-1} \leq M_{1}$;\\
(ii) for some integer $v(k) \in\{k,k-1,\ldots,\max\{k-m+1,0\}\}$,
\begin{equation*}\label{alphack} 
	\alpha_{k} \leq \frac{g_{v(k)}\tr \psi(A) g_{v(k)}}{g_{v(k)}\tr A\psi(A) g_{v(k)}},
\end{equation*}
where $\psi$ is a real analytic function defined as \eqref{fystepbb}.
\begin{theorem}\label{converth}
	Suppose that the sequence $\left\{\left\|g_{k}\right\|\right\}$ is generated by applying the  method \eqref{snew} to solve problem \eqref{eqpro} with  $A=\diag\{\lambda_1,\lambda_2,\ldots,\lambda_n\}$ and $1=\lambda_1<\lambda_2<\ldots<\lambda_n$. Further assume that $|\alpha_{k-1}^{BB1}-\alpha_{k}^{BB1}|>\tilde{\epsilon}$ for some $\tilde{\epsilon}>0$ and $k>2$. Then either $g_{k} = 0$ for some finite $k$ or the sequence  $\left\{\left\|g_{k}\right\|\right\}$ converges to zero $R$-linearly in the sense that
	\begin{equation*}\label{Citheta}
		\lvert g_{k}^{(i)}\rvert \leqslant  c_{i}\mu^{k}, \quad i=1,2, \cdots, n,
	\end{equation*}
	where  ${\mu}=1-\frac{\lambda_{1}}{M_1}$ and
	\begin{equation*}
		\begin{cases} 
			c_{1}=\lvert g_{1}^{(1)}\rvert; \\
			c_{i}=\displaystyle\max\left\{\lvert g_{1}^{(i)}\rvert,\frac{\lvert g_{2}^{(i)}\rvert}{\mu}, \frac{\max \{\nu_{i}, \nu_{i}^{2}\}}{\mu^{2}}\sqrt{\sum_{j=1}^{i-1} c_{j}^{2}}\right\},~~i=2,3, \cdots, n,
		\end{cases}
	\end{equation*}
	with $\nu_{i}=\max \left\{\frac{\lambda_{i}}{\lambda_{1}}-1,1-\frac{\lambda_{i}}{M_1}\right\}$.
\end{theorem}
\begin{proof}
	It is obvious that $\alpha_k\leq\alpha_k^{BB1}$ for all $k>1$. Thus, condition (ii) holds with $v(k)=k-1$, i.e. $m=2$, and $\psi(A)=I$.
	
	Since $u_k, v_k, r_k$ are orthonormal vectors, by Corollary 4.3.16 of \cite{horn2012matrix} we have 
	\begin{equation*}
		\lambda_i\leq\lambda_i({H}_k)\leq\lambda_{i+n-3},~i=1,2,3,
	\end{equation*}
	which gives  $\lambda_1\leq(\alpha^{new}_k)^{-1}\leq\lambda_n$. Notice that $\lambda_1\leq(\alpha_{k}^{BB1})^{-1},(\alpha_{k}^{BB2})^{-1}\leq\lambda_n$. Thus, when $\alpha_{k}=\min\{\alpha_{k-1}^{BB2},\alpha_{k}^{BB2},\alpha_k^{new}\}$ or $\alpha_{k}^{BB1}$, $\alpha_k^{-1}$ falls into $[\lambda_1,\lambda_n]$.

	Now we consider the case $\alpha_{k}=\min\{\alpha_{k-1}^{BB2},\alpha_{k}^{BB2},\alpha_k^{BBQ}\}$.
	It follows from Theorem 2.9 of \cite{hdlbbq} that either
\begin{equation*}
	\big(\alpha_{k}^{BBQ}\big)^{-1}\leq\frac{\phi_2}{\phi_3}=\frac{1}{\alpha_{k}^{BB2}}+\frac{\phi_1}{\phi_3}\alpha_{k}^{BB1}\leq\lambda_n+\frac{2\lambda_n^2}{\tilde{\epsilon}\lambda_1^2},
\end{equation*}
or
\begin{equation*}\label{nbb1lbd2}
	\big(\alpha_{k}^{BBQ}\big)^{-1}\leq\frac{1}{\max\{\alpha_{k}^{BB2},\alpha_{k-1}^{BB2}\}}\leq \lambda_n,
\end{equation*}
which implies that $\alpha_{k}$ satisfies condition (i) with $M_1=\lambda_n+\frac{2\lambda_n^2}{\tilde{\epsilon}\lambda_1^2}$.

Based on the above analysis, we conclude that condition (i) always holds with $M_1$. 
We finish the proof by resorting to Theorem 1 of \cite{li2023note}.
\end{proof}


We compared the method \eqref{snew} with other recent successful gradient methods including:
\begin{itemize}
	\item[(i)] BB \cite{Barzilai1988two}: the original BB method using $\alpha_k^{BB1}$;
	
	\item[(ii)] DY \cite{dai2005analysis}: the Dai--Yuan monotone gradient method;
	
	%
	
	\item[(iii)] ABBmin2 \cite{frassoldati2008new}: a gradient method adaptively using $\alpha_k^{BB1}$ and a short stepsize in \cite{frassoldati2008new};
	
	\item[(iv)] SDC \cite{de2014efficient}: a gradient method takes $h$ SD iterates followed by $s$ steps with the same modified Yuan stepsize in \cite{dai2005analysis};
	
	\item[(v)] BBQ \cite{hdlbbq}: a gradient method adaptively using $\alpha_k^{BB1}$ and the short stepsize  $\min\{\alpha_{k-1}^{BB2},\alpha_{k}^{BB2},\alpha_k^{BBQ}\}$.
	
\end{itemize}
All the compared methods were implemented by Matlab (v.9.0-R2016a) on a PC with an Intel Core i7, 2.9 GHz processor and 8 GB of RAM running Windows 10 system.

We tested the following quadratic problem \cite{yuan2006new}:
\begin{equation}\label{testqp}
	\min_{x\in\mathbb{R}^n}~~f(x)=(x-x^{*})\tr\diag\{v_1,\ldots,v_n\}(x-x^{*}),
\end{equation}
where $x^{*}$ was randomly generated with components in $[-10,10]$ and $v_{j}, j=1,2, \ldots, n$, were generated according to five different distributions listed in Table \ref{tbspe}. 

\begin{table}[th!b]
	\begin{center}
		\caption{Distributions of $v_j$.}\label{tbspe}
		\renewcommand\arraystretch{1.6}
		\footnotesize
		\setlength{\tabcolsep}{0.5ex}
		\begin{tabular}{|c|c|}
			\hline
			\multirow{1}{*}{sets} &\multicolumn{1}{c|}{spectrum} \\
			\hline
			\multirow{1}{*}{1} 
			&$\{v_2,\ldots,v_{n-1}\}\subset(1,\kappa)$,\quad $v_1=1$, $v_n=\kappa$\\
			\hline

			\multirow{2}{*}{2}
			&$v_{j}=1+\left(\kappa- 1\right)s_j$,\quad $s_j\in(0.8,1)$,~~ $j=1, \ldots, \frac{n}{2}$,\\ &$s_j\in(0,0.2)$,~~ $j=\frac{n}{2}+1, \ldots, n$	\\
			\hline
			
			\multirow{1}{*}{3}
			&$\{v_2,\ldots,v_{n/5}\}\subset(1,100)$,\quad $\{v_{n/5+1},\ldots,v_{n-1}\}\subset(\frac{\kappa}{2},\kappa)$,\quad $v_1=1$, $v_n=\kappa$\\
			\hline
			

			\multirow{1}{*}{4}
			&$v_j=\kappa^{\frac{n-j}{n-1}}$\\
			\hline
			
			
			\multirow{1}{*}{5}
			&$\{v_2,\ldots,v_{4n/5}\}\subset(1,100)$,\quad
			$\{v_{4n/5+1},\ldots,v_{n-1}\}\subset(\frac{\kappa}{2},\kappa)$	\\
			\hline
			
		\end{tabular}
	\end{center}
\end{table}

For ABBmin2, BBQ and method \eqref{snew}, we chose $\tau_1$ from $\{0.1,0.2,\ldots,0.9\}$ to get the best performance for each given $\gamma\in\{1,1.02,1.05,1.1,1.2,1.3\}$. For the SDC method, as \cite{de2014efficient}, we compared all combinations of $h$ from $\{8,10,20,30,40,50\}$ and $s$ from $\{2,4,6,8\}$, and selected the pair that leads to the  best performance. We present chosen parameters for the compared methods on different problem sets in Table \ref{tbsparas}.

\begin{table}[htp!b]
	\centering
	\caption{Parameter settings of the compared methods for the five problem sets in Table \ref{tbspe}.}
	\renewcommand\arraystretch{1.1}
	\footnotesize
	\setlength{\tabcolsep}{1.5ex}
	\begin{tabular}{|c|c|c|c|c|}
		\hline
		
		\multicolumn{1}{|c|}{\multirow{2}{*}{problem sets}} 
		& \multicolumn{1}{c|}{\multirow{1}{*}{SDC}} &\multicolumn{1}{c|}{\multirow{1}{*}{ABBmin2}}  &\multicolumn{1}{c|}{\multirow{1}{*}{BBQ}}  &\multicolumn{1}{c|}{\multirow{1}{*}{method \eqref{snew}}}  \\
		
		\cline{2-5}
		&$(h,s)$ &$(\tau_1,\gamma)$  &$(\tau_1,\gamma)$ &$(\tau_1,\gamma)$ \\
		\hline
		
		\multicolumn{1}{|c|}{\multirow{1}{*}{1}} 	&$(50,4)$ & $(0.9,1)$  & $(0.2,1)$  & $(0.9,1)$  \\
		\hline
		
		\multicolumn{1}{|c|}{\multirow{1}{*}{2}} &$(8,8)$ & $(0.9,1)$  & $(0.8,1)$  & $(0.9,1)$   \\
		\hline
		
	\multicolumn{1}{|c|}{\multirow{1}{*}{3}} &$(8,8)$ & $(0.3,1.3)$  & $(0.6,1.3)$  & $(0.5,1)$  \\
\hline

\multicolumn{1}{|c|}{\multirow{1}{*}{4}} &$(50,6)$ & $(0.5,1)$  & $(0.4,1)$  & $(0.5,1)$  \\
\hline

\multicolumn{1}{|c|}{\multirow{1}{*}{5}} &$(8,8)$ & $(0.2,1.3)$  & $(0.3,1.3)$  & $(0.6,1.3)$  \\
\hline		
%
%

	\end{tabular}
	\label{tbsparas}
\end{table}

\begin{table}[ht!b]
	\centering
	\caption{Average number of iterations required by the compared methods on the five problem sets  in Table \ref{tbspe}.}
		\renewcommand\arraystretch{1.3}
		\scriptsize
	\setlength{\tabcolsep}{2.ex}
	\begin{tabular}{|c|c|c|c|c|c|c|c|}
		\hline
		
		\multicolumn{1}{|c|}{\multirow{1}{*}{$\kappa$}} & \multicolumn{1}{c|}{\multirow{1}{*}{$\epsilon$}} & \multicolumn{1}{c|}{\multirow{1}{*}{BB}} & \multicolumn{1}{c|}{\multirow{1}{*}{DY}} & \multicolumn{1}{c|}{\multirow{1}{*}{SDC}} &\multicolumn{1}{c|}{\multirow{1}{*}{ABBmin2}}  &\multicolumn{1}{c|}{\multirow{1}{*}{BBQ}}  &\multicolumn{1}{c|}{\multirow{1}{*}{method \eqref{snew}}}  \\
		
\hline
\multicolumn{8}{|c|}{\multirow{1}{*}{problem set 1}} \\
\hline
\multicolumn{1}{|c|}{\multirow{3}{*}{$10^{4}$}} 
&$10^{-6}$ & 252.9  & 225.2  & 222.7  & 287.4  & 234.2  & 250.8  \\
\cline{2-8}
&$10^{-9}$ & 882.0  & 751.1  & 726.0  & 619.8  & 688.9  & 585.0  \\
\cline{2-8}
&$10^{-12}$ &1492.6  &1212.2  &1212.6  & 943.3  &1046.3  & 879.9  \\
\hline
\multicolumn{1}{|c|}{\multirow{3}{*}{$10^{5}$}} 
&$10^{-6}$ & 230.5  & 193.9  & 194.5  & 194.9  & 195.0  & 198.7  \\
\cline{2-8}
&$10^{-9}$ &1587.2  &1612.8  &1474.3  & 600.4  & 809.1  & 611.1  \\
\cline{2-8}
&$10^{-12}$ &3213.3  &3320.6  &2816.9  & 834.0  &1125.9  & 819.7  \\
\hline
\multicolumn{1}{|c|}{\multirow{3}{*}{$10^{6}$}} 
&$10^{-6}$ & 262.8  & 238.9  & 237.4  & 234.3  & 225.7  & 236.0  \\
\cline{2-8}
&$10^{-9}$ &5187.0  &3861.1  &2367.7  &1021.9  &1579.7  &1151.4  \\
\cline{2-8}
&$10^{-12}$ &12320.1  &12920.7  &5933.6  &1512.2  &1950.8  &1631.9  \\
\hline

\multicolumn{8}{|c|}{\multirow{1}{*}{problem set 2}} \\
\hline
\multicolumn{1}{|c|}{\multirow{3}{*}{$10^{4}$}} 
&$10^{-6}$ & 251.2  & 220.1  & 202.8  & 235.3  & 200.7  & 204.4  \\
\cline{2-8}
&$10^{-9}$ & 686.3  & 575.5  & 493.8  & 504.7  & 464.9  & 470.5  \\
\cline{2-8}
&$10^{-12}$ &1081.7  & 894.7  & 758.9  & 746.2  & 737.6  & 697.2  \\
\hline
\multicolumn{1}{|c|}{\multirow{3}{*}{$10^{5}$}} 
&$10^{-6}$ & 261.1  & 252.6  & 233.8  & 297.7  & 222.0  & 238.5  \\
\cline{2-8}
&$10^{-9}$ &1047.7  & 783.0  & 683.4  & 708.0  & 666.1  & 628.3  \\
\cline{2-8}
&$10^{-12}$ &1974.2  &1413.2  &1085.3  &1090.0  &1020.8  & 958.5  \\
\hline
\multicolumn{1}{|c|}{\multirow{3}{*}{$10^{6}$}} 
&$10^{-6}$ & 284.0  & 266.7  & 250.3  & 314.8  & 240.6  & 245.9  \\
\cline{2-8}
&$10^{-9}$ &1895.5  &1185.8  &1023.7  & 970.4  & 904.8  & 832.3  \\
\cline{2-8}
&$10^{-12}$ &3369.9  &2267.5  &1619.2  &1525.8  &1387.2  &1295.5  \\
\hline

\multicolumn{8}{|c|}{\multirow{1}{*}{problem set 3}} \\
\hline
\multicolumn{1}{|c|}{\multirow{3}{*}{$10^{4}$}} 
&$10^{-6}$ & 257.0  & 198.8  & 203.9  & 147.2  & 157.3  & 140.5  \\
\cline{2-8}
&$10^{-9}$ & 613.3  & 525.5  & 506.5  & 372.1  & 356.3  & 337.3  \\
\cline{2-8}
&$10^{-12}$ & 953.3  & 864.1  & 731.6  & 555.5  & 532.8  & 527.3  \\
\hline
\multicolumn{1}{|c|}{\multirow{3}{*}{$10^{5}$}} 
&$10^{-6}$ & 351.4  & 321.3  & 166.5  &  93.6  &  94.6  &  96.3  \\
\cline{2-8}
&$10^{-9}$ &1291.0  &1233.4  & 730.8  & 388.2  & 382.8  & 366.0  \\
\cline{2-8}
&$10^{-12}$ &2446.1  &2137.4  &1263.9  & 644.6  & 619.8  & 633.5  \\
\hline
\multicolumn{1}{|c|}{\multirow{3}{*}{$10^{6}$}} 
&$10^{-6}$ & 400.3  & 259.0  &  93.5  &  48.9  &  50.5  &  49.5  \\
\cline{2-8}
&$10^{-9}$ &2674.6  &2203.8  & 863.3  & 355.8  & 372.2  & 357.6  \\
\cline{2-8}
&$10^{-12}$ &5931.4  &4487.3  &1456.0  & 668.4  & 655.8  & 619.7  \\
\hline

\multicolumn{8}{|c|}{\multirow{1}{*}{problem set 4}} \\
\hline
\multicolumn{1}{|c|}{\multirow{3}{*}{$10^{4}$}} 
&$10^{-6}$ & 609.7  & 563.9  & 511.1  & 481.5  & 462.6  & 469.6  \\
\cline{2-8}
&$10^{-9}$ &1239.5  & 999.9  & 971.1  & 877.8  & 892.1  & 849.9  \\
\cline{2-8}
&$10^{-12}$ &1771.3  &1387.0  &1409.8  &1262.7  &1278.9  &1257.6  \\
\hline
\multicolumn{1}{|c|}{\multirow{3}{*}{$10^{5}$}} 
&$10^{-6}$ &1588.9  &1237.1  &1158.6  &1133.2  &1124.2  &1042.1  \\
\cline{2-8}
&$10^{-9}$ &3619.5  &2895.9  &2826.6  &2412.1  &2428.8  &2361.9  \\
\cline{2-8}
&$10^{-12}$ &5655.2  &4640.3  &4283.5  &3623.0  &3711.8  &3581.6  \\
\hline
\multicolumn{1}{|c|}{\multirow{3}{*}{$10^{6}$}} 
&$10^{-6}$ &2578.4  &2348.0  &2095.7  &2142.4  &2092.5  &1935.1  \\
\cline{2-8}
&$10^{-9}$ &9873.7  &11463.4  &7541.1  &6191.6  &6634.6  &6210.3  \\
\cline{2-8}
&$10^{-12}$ &18552.5  &20106.6  &11403.3  &10162.9  &10456.4  &9970.3  \\
\hline

\multicolumn{8}{|c|}{\multirow{1}{*}{problem set 5}} \\
\hline
\multicolumn{1}{|c|}{\multirow{3}{*}{$10^{4}$}} 
&$10^{-6}$ & 313.8  & 299.7  & 263.3  & 186.3  & 196.8  & 182.5  \\
\cline{2-8}
&$10^{-9}$ & 692.8  & 648.9  & 511.5  & 383.0  & 386.0  & 380.7  \\
\cline{2-8}
&$10^{-12}$ & 943.9  & 972.9  & 741.4  & 591.0  & 577.0  & 555.2  \\
\hline
\multicolumn{1}{|c|}{\multirow{3}{*}{$10^{5}$}} 
&$10^{-6}$ & 505.3  & 504.9  & 275.9  & 172.3  & 153.9  & 152.1  \\
\cline{2-8}
&$10^{-9}$ &1637.6  &1463.7  & 829.6  & 444.1  & 442.0  & 409.3  \\
\cline{2-8}
&$10^{-12}$ &2580.5  &2352.3  &1318.7  & 681.3  & 656.9  & 658.2  \\
\hline
\multicolumn{1}{|c|}{\multirow{3}{*}{$10^{6}$}} 
&$10^{-6}$ & 687.3  & 449.7  & 185.8  &  79.6  &  78.4  &  79.7  \\
\cline{2-8}
&$10^{-9}$ &3431.9  &2729.0  & 972.0  & 411.1  & 388.1  & 388.0  \\
\cline{2-8}
&$10^{-12}$ &6012.9  &5090.1  &1877.4  & 732.7  & 708.1  & 668.4  \\
\hline

	\end{tabular}
	\label{tbs2}
\end{table}

The problem dimension was set to $n=10000$ in this test. For each method, the iteration was stopped once the gradient norm reduces by a factor of $\epsilon$ or exceeds the maximum iteration number 50000. We set $\kappa$ to $10^4, 10^5, 10^6$ and $\epsilon$ to $10^{-6}, 10^{-9}, 10^{-12}$ to see performances of the compared methods on different values of  condition number and tolerance. For each value of $\kappa$ or $\epsilon$, average results over 10 different starting points with entries randomly generated in $[-10,10]$ are presented in Table \ref{tbs2}. 


We can see that for most of the problem sets our method \eqref{snew} performs much better than the BB, DY and SDC methods, and better than the ABBmin2 and BBQ methods in the sense of number of iterations. It is worth to mention that the DY, SDC, ABBmin2 and BBQ methods employ stepsizes with the two-dimensional quadratic termination property. Moreover, applying the DY, SDC and ABBmin2 methods to general unconstrained optimization is not easy since computing their stepsizes requires either exact line searches or the Hessian. 


\subsection{Unconstrained optimization}\label{secunc}
In this subsection, we extend the method \eqref{snew} to solve general unconstrained optimization 
 \begin{equation}\label{equncpro}
	 	\min_{x\in\mathbb{R}^n}~~f(x),
	 \end{equation}
 where $f: \mathbb{R}^n \to \mathbb{R}$ is continuously differentiable.

Line searches are essential to ensure global convergence of the BB method since directly apply it to the unconstrained optimization problem \eqref{equncpro} may not converge even though the objective function is strongly convex, see \cite{burdakov2019stabilized} for example. Due to the inherent nonmonotone property of BB-like stepsizes, nonmonotone line searches \cite{dai2005projected,dai2001adaptive,grippo1986nonmonotone,zhang2004nonmonotone} are often employed to achieve good performance.  
%
Here we adopt the Dai--Fletcher nonmonotone line search \cite{dai2005projected} which
finds a stepsize such that
\begin{equation}\label{DFcond}
	f(x_k+\lambda_k d_k)\leq f_{r}+\delta \lambda_k g_k\tr d_k, 
\end{equation}
where  $\delta\in(0,1)$ and $f_r$ is the reference value. 


To update the reference function value $f_r$, the Dai--Fletcher nonmonotone line search needs the current best function value $f_{\min}=\min_{1\leq j\leq k} f(x_{j})$, the maximal value of the objective function since $f_{\min}$ was found, namely $f_{c}$, and a preassigned number $T$. Denote $t$ by the number of iterations since the value of $f_{\min}$ was actually obtained. We initialize $f_{r}=f_{\min}=f_c=f(x_1)$, and $t=0$. If the method can find a better function value in $T$ iterations, where $T$ is a preset positive integer constant, then the value of $f_r$ remains unchanged. Otherwise, if $t$ reaches $T$, the reference function value $f_r$ is reset to a candidate value $f_c$ and $f_c$ is reset to the current function value.  
The strategy to update $f_r$ can be summarized in Matlab-style notation as follows.
\begin{align*}
	&\mbox{if} \ f_k< f_{\min} \\
	&\quad\quad f_{\min}=f_k, \ \ f_c=f_k, \ \ t=0, \\
	&\mbox{else} \\
	&\quad\quad f_c=\max\{f_c,f_k\}, \ \ t=t+1, \\
	&\quad\quad \textsf{if} \ t=T\\
	&\quad \quad\quad f_{r}=f_c,\ \ f_c=f_k, \ t=0,\\
	&\quad\quad \mbox{end}\\
	&\ \mbox{end}
\end{align*}

Our method for unconstrained optimization \eqref{equncpro} is presented in Algorithm \ref{alunc}. Under standard assumptions, we can establish global convergence  of Algorithm \ref{alunc} in the sense $\lim\inf_{k\rightarrow\infty}\|g_k\|=0$ and $R$-linear rate for strongly convex objective functions, see Theorem 3.2 in \cite{dai2001adaptive} and Theorem 4.1 in \cite{huang2015rate}.

\begin{algorithm}[htp!b] 
	\caption{A gradient method for unconstrained optimization} 
	\label{alunc}
		\begin{algorithmic} 
		\STATE{Initialization: $x_{1}\in \mathbb{R}^n$, $\alpha_{\max}\geq\alpha_{\min}>0$, $\alpha_{1}\in[\alpha_{\min},\alpha_{\max}]$,  $\tau_{1}>0$, $\gamma\geq1$, $\epsilon, \delta, \eta\in(0,1)$, $T\in\mathbb{N}$, $k:=1$.}
		 
		\WHILE{$\left\|g_{k}\right\|>\epsilon$}
			\STATE{$d_{k}=-g_{k}$, $\lambda_{k}=\alpha_k$}

			\WHILE{the condition \eqref{DFcond} does not hold}
				\STATE{$\lambda_{k}=\eta\lambda_{k}$}
			\ENDWHILE
			
			\STATE{$x_{k+1}=x_{k}+\lambda_{k}d_{k}$}
			
			\IF{$s_{k}\tr  y_{k}>0$}
				\IF{$\alpha_k^{BB2}/\alpha_k^{BB1}<\tau_k$}		
					\IF{$s_{k-1}\tr  y_{k-1}>0$ and $s_{k-2}\tr  y_{k-2}>0$}
						\IF{$\alpha_{k+1}^{new}>0$}	
							\STATE{${\alpha}_{k+1}=\min\{\alpha_{k}^{BB2}, \alpha_{k+1}^{BB2}, \alpha_{k+1}^{new}\}$}
						\ELSE
							\STATE{$\alpha_{k+1} =\min\{\alpha_{k}^{BB2}, \alpha_{k+1}^{BB2}\}$}
						\ENDIF			
					\ELSIF{$s_{k-1}\tr  y_{k-1}>0$}
						\IF{$\alpha_{k+1}^{BBQ}>0$}
							\STATE{${\alpha}_{k+1}=\min\{\alpha_{k}^{BB2}, \alpha_{k+1}^{BB2}, \alpha_{k+1}^{BBQ}\}$}	
						\ELSE
							\STATE{$\alpha_{k+1} =\min\{\alpha_{k}^{BB2}, \alpha_{k+1}^{BB2}\}$}
						\ENDIF
					\ELSE
						\STATE{$\alpha_{k+1} =\alpha_{k+1}^{BB2}$}
					\ENDIF	
					\STATE{$\tau_{k+1}=\tau_k/\gamma$}
			 \ELSE
				\STATE{$\alpha_{k+1}=\alpha_{k+1}^{BB1}$}
				\STATE{$\tau_{k+1}=\tau_k\gamma$}
			\ENDIF
		\ELSE
			\STATE{$\alpha_{k+1} = \min\{1/\|g_{k+1}\|_\infty,\|x_{k+1}\|_\infty/\|g_{k+1}\|_\infty\}$}
		\ENDIF
			
		\STATE{	Chop extreme values of the stepsize such that $\alpha_{k+1}\in[\alpha_{\min},\alpha_{\max}]$}
			
		\STATE{Compute $\tau_{k+1}$ by \eqref{tauk}}
			
		\STATE{$k:=k+1$}
		\ENDWHILE
	\end{algorithmic} 
\end{algorithm}

We report numerical results of Algorithm \ref{alunc} on general unconstrained problems from the CUTEst collection \cite{gould2015cutest} with dimension less than or equal to $10000$.  As shown in \cite{di2018steplength,hdlbbq}, the ABBmin method \cite{di2018steplength,frassoldati2008new} is very competitive with the LMSD method \cite{fletcher2012limited} while the BBQ method is much faster than the GBB \cite{raydan1997barzilai} and ABBmin methods. Hence we only compared  Algorithm \ref{alunc} with the BBQ method. To this end, we built a virtual machine with 4 GB of RAM running Ubuntu 16 system on the PC mentioned above and implemented the methods by Matlab (v.9.3-R2017b).

We deleted the problem if either it can not be solved in $200000$ iterations by any of the algorithms or the function evaluation exceeds one million, and $162$ problems are left.

\begin{figure}[htp!b]
	\centering
	\subfloat{\label{fig:2a}\includegraphics[width=0.47\textwidth,height=0.35\textwidth]{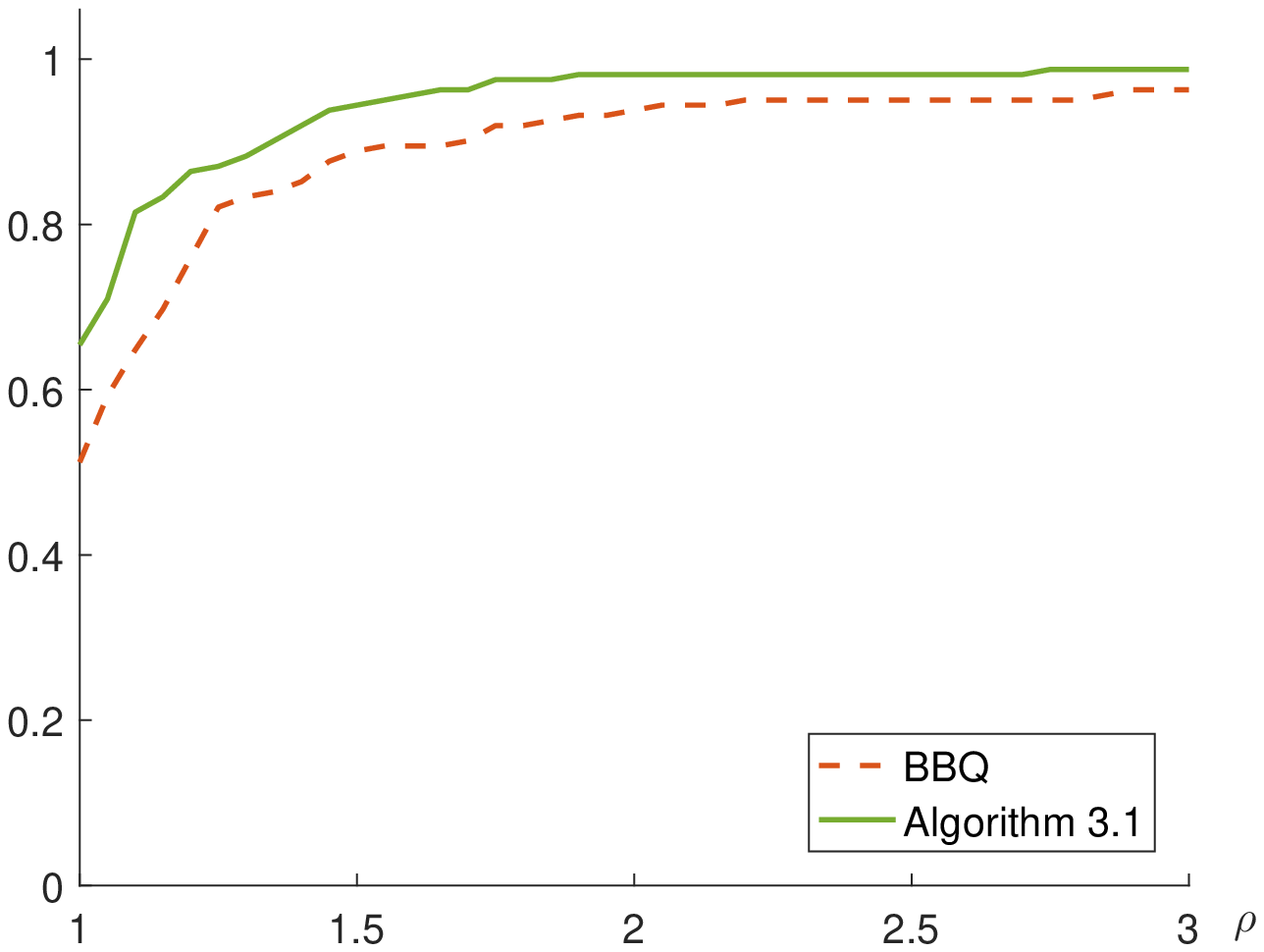}}
	\subfloat{\label{fig:2b}\includegraphics[width=0.47\textwidth,height=0.35\textwidth]{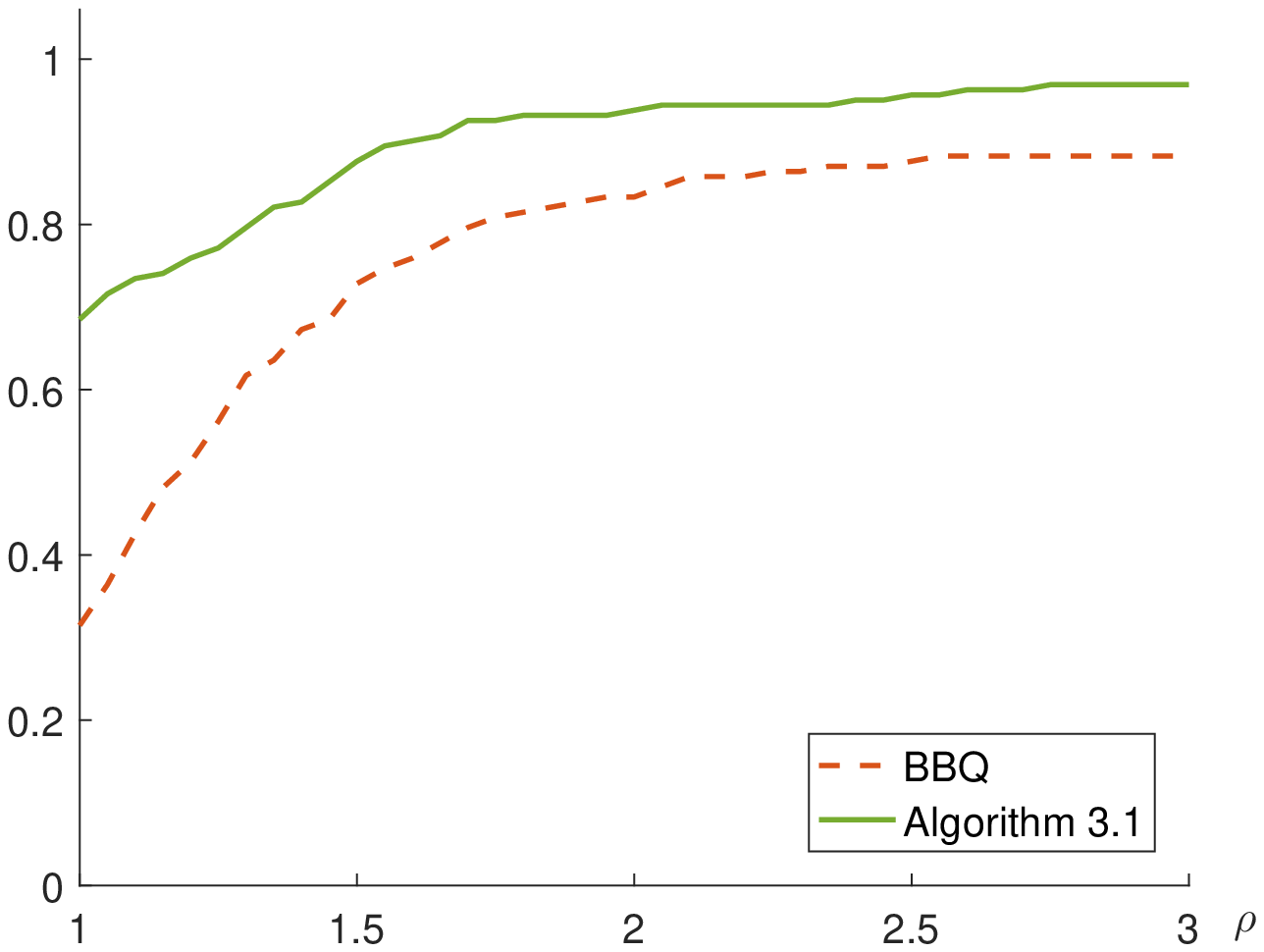}}
	\caption{Performance profiles of Algorithm \ref{alunc} and BBQ on 162 unconstrained problems from CUTEst, iteration (left) and CPU time (right) metrics.}\label{cutestunc}
\end{figure}

For Algorithm \ref{alunc}, similar to the BBQ method \cite{hdlbbq}, we used $\alpha_{1}=\|x_{1}\|_\infty/\|g_{1}\|_\infty$ if $x_{1}\neq0$ and otherwise $\alpha_{1}=1/\|g_{1}\|_\infty$. In addition, we set $\alpha_{i}=\alpha_{i}^{BB1}$ for the case $\alpha_{i}^{BB1}>0$ and otherwise chose $\alpha_{i} = \min\{1/\|g_{i}\|_\infty,\|x_{i}\|_\infty/\|g_{i}\|_\infty\}$, $i=2,3,4$. The following parameters were employed  for Algorithm \ref{alunc}:
$$\alpha_{\min}=10^{-10},\, \alpha_{\max}=10^{6},\, T=3,\, \delta=10^{-4},\,\eta=0.5,\, \tau_1=0.65, \, \gamma=1.4.$$
Default parameter settings were used for the BBQ method. We stopped the iteration of the two methods if $\|g_{k}\|_\infty\leq10^{-6}$.

Performance profiles \cite{dolan2002} of Algorithm \ref{alunc} and the BBQ method using iteration and CPU time metrics are plotted in Figure \ref{cutestunc}. For each method, the vertical axis of the figure shows the percentage of problems the method solves within the factor $\rho$ of the minimum value of the metric. From Figure \ref{cutestunc}, it can be seen that Algorithm \ref{alunc} performs better than the BBQ method.

\begin{table}[htp!b]
	\setlength{\tabcolsep}{1.1ex}
	\caption{Results of Algorithm \ref{alunc} and BBQ on 162 unconstrained problems from the CUTEst collection.}\label{tbunc}
	\centering
	\begin{scriptsize}
		\begin{tabular}{|c|c|cccc|cccc|}
			\hline
			\multicolumn{1}{|c|}{\multirow{2}{*}{problem}} &\multicolumn{1}{c|}{\multirow{2}{*}{$n$}} &\multicolumn{4}{c|}{\multirow{1}{*}{BBQ}} &\multicolumn{4}{c|}{\multirow{1}{*}{Algorithm \ref{alunc}}}\\
			\cline{3-10}
			&	&nfe &ngrad  &iter &time &nfe &ngrad &iter &time\\		
			\hline
			AKIVA   &     2	&    34	&    25	&    25	&0.4641   &    31	 &    24	&    24  &0.0744 \\
			ALLINITU   &     4	&    16	&    16	&    16	&0.0387   &    16	 &    16	&    16  &0.0239 \\
			ARGLINA   &   200	&     3	&     3	&     3	&0.0063   &     3	 &     3	&     3  &0.0019 \\
			ARGTRIGLS   &   200	&  1809	&  1581	&  1581	&1.5536   &  2358	 &  2104	&  2104  &1.4290 \\
			ARWHEAD   &  5000	&     4	&     4	&     4	&0.0182   &     4	 &     4	&     4  &0.0361 \\
			BA-L1LS   &    57	&    43	&    32	&    32	&0.0027   &    45	 &    34	&    34  &0.0172 \\
			BARD   &     3	&    85	&    73	&    73	&0.0072   &    84	 &    80	&    80  &0.0129 \\
			BDQRTIC   &  5000	&    74	&    73	&    73	&0.1262   &    68	 &    68	&    68  &0.1069 \\
			BEALE   &     2	&    40	&    40	&    40	&0.0026   &    32	 &    32	&    32  &0.0023 \\
			BENNETT5LS   &     3	&   158	&   138	&   138	&0.0396   &    56	 &    44	&    44  &0.0123 \\
			BIGGS6   &     6	&   407	&   343	&   343	&0.0386   &   387	 &   361	&   361  &0.0239 \\
			BOX   & 10000	&    53	&    26	&    26	&0.1999   &    50	 &    24	&    24  &0.2215 \\
			BOX3   &     3	&    38	&    31	&    31	&0.0037   &    34	 &    32	&    32  &0.0035 \\
			BOXBODLS   &     2	&     9	&     8	&     8	&0.0060   &    73	 &    72	&    72  &0.0140 \\
			BRKMCC   &     2	&    10	&     9	&     9	&0.0010   &    10	 &     9	&     9  &0.0008 \\
			BROWNAL   &   200	&   602	&   366	&   366	&0.1763   &    64	 &    63	&    63  &0.0237 \\
			BROWNDEN   &     4	&    69	&    65	&    65	&0.0038   &    66	 &    63	&    63  &0.0050 \\
			BROYDN3DLS   &  5000	&    90	&    88	&    88	&0.0843   &    92	 &    92	&    92  &0.1298 \\
			BROYDN7D   &  5000	&  3264	&  2592	&  2592	&11.0259   &  2690	 &  2461	&  2461  &10.1276 \\
			BROYDNBDLS   &  5000	&    50	&    49	&    49	&0.3051   &    41	 &    40	&    40  &0.0900 \\
			BRYBND   &  5000	&    50	&    49	&    49	&0.1402   &    41	 &    40	&    40  &0.0955 \\
			CHAINWOO   &  4000	&   554	&   492	&   492	&0.7618   &   550	 &   522	&   522  &0.7527 \\
			CHNROSNB   &    50	&  1032	&   869	&   869	&0.0781   &   707	 &   637	&   637  &0.0461 \\
			CHNRSNBM   &    50	&  1264	&  1068	&  1068	&0.0698   &   832	 &   750	&   750  &0.0452 \\
			CHWIRUT1LS   &     3	&   143	&   112	&   112	&0.0185   &    95	 &    89	&    89  &0.0126 \\
			CHWIRUT2LS   &     3	&    86	&    67	&    67	&0.0043   &   107	 &    91	&    91  &0.0072 \\
			COSINE   & 10000	&    22	&    22	&    22	&0.2445   &    24	 &    24	&    24  &0.1507 \\
			CRAGGLVY   &  5000	&   404	&   341	&   341	&0.9754   &   160	 &   157	&   157  &0.4774 \\
			CUBE   &     2	&   152	&   136	&   136	&0.0199   &    47	 &    47	&    47  &0.0052 \\
			DECONVU   &    63	&  4562	&  3983	&  3983	&0.3496   &  3747	 &  3279	&  3279  &0.3670 \\
			DENSCHNA   &     2	&    15	&    15	&    15	&0.0019   &    13	 &    13	&    13  &0.0013 \\
			DENSCHNB   &     2	&     9	&     9	&     9	&0.0012   &     9	 &     9	&     9  &0.0009 \\
			DENSCHNC   &     2	&     7	&     7	&     7	&0.0006   &     7	 &     7	&     7  &0.0005 \\
			DENSCHND   &     3	&   128	&   121	&   121	&0.0090   &   112	 &   112	&   112  &0.0060 \\
			DENSCHNE   &     3	&    71	&    49	&    49	&0.0048   &    41	 &    40	&    40  &0.0038 \\
			DENSCHNF   &     2	&    18	&    18	&    18	&0.0019   &    15	 &    15	&    15  &0.0015 \\
			DIXMAANA   &  3000	&     8	&     8	&     8	&0.0262   &     8	 &     8	&     8  &0.0232 \\
			DIXMAANB   &  3000	&     8	&     8	&     8	&0.0207   &     8	 &     8	&     8  &0.0260 \\
			DIXMAANC   &  3000	&     9	&     9	&     9	&0.0279   &     9	 &     9	&     9  &0.0080 \\
			DIXMAAND   &  3000	&    10	&    10	&    10	&0.0299   &     9	 &     9	&     9  &0.0098 \\
			DIXMAANE   &  3000	&   417	&   373	&   373	&0.3061   &   350	 &   331	&   331  &0.2373 \\
			DIXMAANF   &  3000	&   346	&   305	&   305	&0.2789   &   342	 &   329	&   329  &0.2412 \\
			DIXMAANG   &  3000	&   423	&   375	&   375	&0.3209   &   271	 &   257	&   257  &0.2629 \\
			DIXMAANH   &  3000	&   314	&   281	&   281	&0.2967   &   358	 &   335	&   335  &0.2691 \\
			DIXMAANI   &  3000	&  2487	&  2143	&  2143	&1.6983   &  3573	 &  3144	&  3144  &2.2829 \\
			DIXMAANJ   &  3000	&   355	&   318	&   318	&0.2676   &   327	 &   314	&   314  &0.2354 \\
			DIXMAANK   &  3000	&   314	&   283	&   283	&0.2349   &   306	 &   294	&   294  &0.3479 \\
			DIXMAANL   &  3000	&   260	&   238	&   238	&0.2277   &   210	 &   200	&   200  &0.1660 \\
			DIXMAANM   &  3000	&  3267	&  2842	&  2842	&2.0558   &  2915	 &  2561	&  2561  &1.9314 \\
			DIXMAANN   &  3000	&   788	&   684	&   684	&0.6266   &   900	 &   794	&   794  &0.5834 \\
			DIXMAANO   &  3000	&   858	&   759	&   759	&0.6014   &   670	 &   612	&   612  &0.4563 \\
			DIXMAANP   &  3000	&   645	&   575	&   575	&0.9761   &   916	 &   828	&   828  &0.6306 \\
			DIXON3DQ   & 10000	&  6062	&  5285	&  5285	&14.6310   &  5782	 &  5126	&  5126  &21.1504 \\
			DJTL   &     2	&   393	&   191	&   191	&0.0183   &   395	 &   191	&   191  &0.0188 \\
			DQDRTIC   &  5000	&    15	&    15	&    15	&0.0145   &    17	 &    17	&    17  &0.0228 \\
			DQRTIC   &  5000	&    50	&    50	&    50	&0.0348   &    53	 &    53	&    53  &0.0360 \\
			ECKERLE4LS   &     3	&    13	&     5	&     5	&0.0007   &    13	 &     5	&     5  &0.0017 \\
			EDENSCH   &  2000	&    38	&    38	&    38	&0.0199   &    41	 &    41	&    41  &0.0154 \\
			EG2   &  1000	&     6	&     6	&     6	&0.0027   &     6	 &     6	&     6  &0.0024 \\
			EIGENALS   &  2550	& 37898	& 33151	& 33151	&262.6693   & 28858	 & 25585	& 25585  &215.6744 \\
			ENGVAL1   &  5000	&    36	&    30	&    30	&0.0526   &    37	 &    31	&    31  &0.0720 \\
			ENGVAL2   &     3	&   106	&    93	&    93	&0.0065   &   108	 &   107	&   107  &0.0061 \\
			ENSOLS   &     9	&    98	&    88	&    88	&0.0198   &    94	 &    89	&    89  &0.0296 \\
			ERRINROS   &    50	&  2942	&  2605	&  2605	&0.1259   &  3607	 &  3212	&  3212  &0.1637 \\
			ERRINRSM   &    50	& 12125	&  9111	&  9111	&0.4304   &  7243	 &  6303	&  6303  &0.2744 \\
			EXPFIT   &     2	&    38	&    34	&    34	&0.0096   &    33	 &    32	&    32  &0.0026 \\
			EXTROSNB   &  1000	& 28845	& 26218	& 26218	&2.9969   & 24148	 & 19874	& 19874  &2.3297 \\
				FBRAINLS   &     2	&    27	&    20	&    20	&0.0325   &    30	 &    23	&    23  &0.0495 \\
			FLETBV3M   &  5000	&    59	&    58	&    58	&0.2234   &    88	 &    88	&    88  &0.6079 \\
		
			\hline
		\end{tabular}
	\end{scriptsize}
\end{table}

\begin{table}[htp!b]
	\setlength{\tabcolsep}{0.8ex}
	\caption{Results of Algorithm \ref{alunc} and BBQ on 162 unconstrained problems from the CUTEst collection (continued).}\label{tbunc2}
	\centering
	\begin{scriptsize}
		\begin{tabular}{|c|c|cccc|cccc|}
			\hline
			\multicolumn{1}{|c|}{\multirow{2}{*}{problem}} &\multicolumn{1}{c|}{\multirow{2}{*}{$n$}} &\multicolumn{4}{c|}{\multirow{1}{*}{BBQ}} &\multicolumn{4}{c|}{\multirow{1}{*}{Algorithm \ref{alunc}}}\\
			\cline{3-10}
			&	&nfe &ngrad  &iter &time &nfe &ngrad &iter &time\\		
			\hline
			
				FLETCBV2   &  5000	&     1	&     1	&     0	&0.0027   &     1	 &     1	&     0  &0.0039 \\
			FLETCHCR   &  1000	&   299	&   270	&   270	&0.0408   &   258	 &   238	&   238  &0.0404 \\
			FMINSRF2   &  5625	&   681	&   622	&   622	&1.3331   &   542	 &   506	&   506  &1.4096 \\
			FMINSURF   &  5625	&   827	&   766	&   766	&2.0244   &   710	 &   655	&   655  &1.4517 \\
			FREUROTH   &  5000	&    62	&    55	&    55	&0.1442   &    63	 &    60	&    60  &0.1833 \\
			GAUSS2LS   &     8	&  2589	&  1334	&  1334	&0.2539   &  3260	 &  1388	&  1388  &0.2710 \\
			GBRAINLS   &     2	&    25	&    17	&    17	&0.0564   &    31	 &    22	&    22  &0.0553 \\
			GENHUMPS   &  5000	& 11216	&  8981	&  8981	&73.2407   &  6265	 &  6005	&  6005  &42.4806 \\
			GENROSE   &   500	&  3962	&  3074	&  3074	&0.3494   &  3287	 &  2822	&  2822  &0.2815 \\
			GROWTHLS   &     3	&     2	&     2	&     2	&0.0011   &     2	 &     2	&     2  &0.0002 \\
			GULF   &     3	&   442	&   346	&   346	&0.0564   &   280	 &   244	&   244  &0.0388 \\
			HAIRY   &     2	&    48	&    44	&    44	&0.0019   &    84	 &    75	&    75  &0.0207 \\
			HATFLDD   &     3	&    50	&    37	&    37	&0.0018   &   102	 &   101	&   101  &0.0083 \\
			HATFLDE   &     3	&   214	&   151	&   151	&0.0175   &   129	 &   127	&   127  &0.0085 \\
			HATFLDFL   &     3	&   203	&   169	&   169	&0.0063   &    67	 &    60	&    60  &0.0040 \\
			HEART8LS   &     8	&   959	&   699	&   699	&0.0462   &   629	 &   587	&   587  &0.0298 \\
			HELIX   &     3	&    69	&    64	&    64	&0.0031   &    38	 &    38	&    38  &0.0022 \\
			HILBERTA   &     2	&     8	&     8	&     8	&0.0008   &     8	 &     8	&     8  &0.0004 \\
			HILBERTB   &    10	&     8	&     8	&     8	&0.0006   &     8	 &     8	&     8  &0.0081 \\
			HIMMELBB   &     2	&     2	&     2	&     2	&0.0005   &     2	 &     2	&     2  &0.0001 \\
			HIMMELBF   &     4	&   363	&   317	&   317	&0.0160   &   549	 &   445	&   445  &0.0266 \\
			HIMMELBG   &     2	&    12	&    11	&    11	&0.0009   &    12	 &    11	&    11  &0.0005 \\
			HIMMELBH   &     2	&    14	&    13	&    13	&0.0013   &    12	 &    11	&    11  &0.0009 \\
			HUMPS   &     2	&   255	&   234	&   234	&0.0141   &   166	 &   166	&   166  &0.0069 \\
			INTEQNELS   &    12	&     9	&     8	&     8	&0.0006   &     9	 &     8	&     8  &0.0005 \\
			JENSMP   &     2	&    31	&    28	&    28	&0.0015   &    53	 &    48	&    48  &0.0020 \\
			KOWOSB   &     4	&   178	&   139	&   139	&0.0101   &   132	 &   122	&   122  &0.0082 \\
			LANCZOS1LS   &     6	&  3148	&  2923	&  2923	&0.1349   &  1623	 &  1488	&  1488  &0.0795 \\
			LANCZOS2LS   &     6	&  2810	&  2613	&  2613	&0.1238   &  2691	 &  2407	&  2407  &0.1269 \\
			LANCZOS3LS   &     6	&  2492	&  2321	&  2321	&0.1104   &  1511	 &  1336	&  1336  &0.0795 \\
			LIARWHD   &  5000	&    51	&    48	&    48	&0.0667   &    47	 &    46	&    46  &0.0960 \\
			LOGHAIRY   &     2	&   260	&   247	&   247	&0.0149   &   411	 &   404	&   404  &0.0240 \\
			LSC1LS   &     3	&   849	&   580	&   580	&0.0261   &   584	 &   579	&   579  &0.0372 \\
			LSC2LS   &     3	&203913	&162654	&162654	&5.7069   &182570	 &132454	&132454  &5.2320 \\
			LUKSAN11LS   &   100	&  5193	&  4173	&  4173	&0.2291   &  3946	 &  3446	&  3446  &0.1748 \\
			LUKSAN12LS   &    98	&   867	&   755	&   755	&0.0787   &   803	 &   635	&   635  &0.0437 \\
			LUKSAN13LS   &    98	&   247	&   192	&   192	&0.0134   &   175	 &   169	&   169  &0.0091 \\
			LUKSAN14LS   &    98	&   173	&   154	&   154	&0.0118   &   207	 &   201	&   201  &0.0306 \\
			LUKSAN15LS   &   100	&    34	&    32	&    32	&0.0251   &    31	 &    29	&    29  &0.0170 \\
			LUKSAN16LS   &   100	&    41	&    38	&    38	&0.0040   &    41	 &    38	&    38  &0.0080 \\
			LUKSAN17LS   &   100	&   494	&   444	&   444	&0.1349   &   469	 &   461	&   461  &0.1208 \\
			LUKSAN21LS   &   100	&  2041	&  1651	&  1651	&0.0908   &  1572	 &  1393	&  1393  &0.0653 \\
			LUKSAN22LS   &   100	&  3843	&  3420	&  3420	&0.2710   &  4165	 &  3609	&  3609  &0.2636 \\
			MANCINO   &   100	&    15	&    15	&    15	&0.0660   &    12	 &    12	&    12  &0.0533 \\
			MARATOSB   &     2	&  3049	&  1487	&  1487	&0.0680   &   220	 &   170	&   170  &0.0126 \\
			MEXHAT   &     2	&    49	&    37	&    37	&0.0048   &    57	 &    50	&    50  &0.0047 \\
			MGH09LS   &     4	&   103	&    92	&    92	&0.0083   &    61	 &    61	&    61  &0.0072 \\
			MNISTS0LS   &   494	&     2	&     2	&     2	&0.5330   &     2	 &     2	&     2  &0.2167 \\
			MNISTS5LS   &   494	&     2	&     2	&     2	&1.0041   &     2	 &     2	&     2  &0.2999 \\
			MOREBV   &  5000	&   119	&    95	&    95	&0.2482   &   108	 &    90	&    90  &0.1706 \\
			MSQRTALS   &  1024	&  5042	&  4468	&  4468	&5.7136   &  4587	 &  4184	&  4184  &5.1710 \\
			MSQRTBLS   &  1024	&  3006	&  2675	&  2675	&3.3117   &  3824	 &  3410	&  3410  &4.9177 \\
			NCB20   &  5010	&   485	&   441	&   441	&2.1745   &   690	 &   638	&   638  &3.0800 \\
			NCB20B   &  5000	&  3132	&  2742	&  2742	&13.3211   &  3979	 &  3150	&  3150  &16.0990 \\
			NONCVXU2   &  5000	& 20662	& 17900	& 17900	&62.0413   & 20667	 & 18286	& 18286  &62.3537 \\
			NONDIA   &  5000	&    25	&    25	&    25	&0.0342   &    39	 &    39	&    39  &0.0522 \\
			NONDQUAR   &  5000	&  3140	&  2706	&  2706	&2.5962   &  2997	 &  2616	&  2616  &2.4613 \\
			OSBORNEA   &     5	& 12518	& 10022	& 10022	&0.4810   &  6795	 &  5854	&  5854  &0.2870 \\
			PALMER5C   &     6	&    28	&    28	&    28	&0.0012   &    20	 &    20	&    20  &0.0009 \\
				PENALTY1   &  1000	&    45	&    45	&    45	&0.0112   &    44	 &    44	&    44  &0.0093 \\
			POWELLSG   &  5000	&   288	&   230	&   230	&0.2233   &   133	 &   125	&   125  &0.1094 \\
			POWER   & 10000	&   891	&   691	&   691	&1.3307   &   872	 &   731	&   731  &1.6106 \\
			QUARTC   &  5000	&    50	&    50	&    50	&0.0505   &    53	 &    53	&    53  &0.0427 \\
			RAT42LS   &     3	&     2	&     2	&     2	&0.0052   &     2	 &     2	&     2  &0.0009 \\
			RAT43LS   &     4	&     4	&     4	&     4	&0.0012   &     4	 &     4	&     4  &0.0004 \\
			ROSENBR   &     2	&   110	&    97	&    97	&0.0055   &    59	 &    57	&    57  &0.0049 \\
				ROSENBRTU   &     2	&   192	&   125	&   125	&0.0097   &   187	 &   167	&   167  &0.0092 \\
				\hline
			
		\end{tabular}		
	\end{scriptsize}	
\end{table}

\begin{table}[htp!b]
	\setlength{\tabcolsep}{0.8ex}
	\caption{Results of Algorithm \ref{alunc} and BBQ on 162 unconstrained problems from the CUTEst collection (continued).}\label{tbunc3}
	\centering
	\begin{scriptsize}
		\begin{tabular}{|c|c|cccc|cccc|}
			\hline
			\multicolumn{1}{|c|}{\multirow{2}{*}{problem}} &\multicolumn{1}{c|}{\multirow{2}{*}{$n$}} &\multicolumn{4}{c|}{\multirow{1}{*}{BBQ}} &\multicolumn{4}{c|}{\multirow{1}{*}{Algorithm \ref{alunc}}}\\
			\cline{3-10}
			&	&nfe &ngrad  &iter &time &nfe &ngrad &iter &time\\		
			\hline

			S308   &     2	&    14	&    14	&    14	&0.0013   &    15	 &    15	&    15  &0.0010 \\
			SCHMVETT   &  5000	&    67	&    67	&    67	&0.2245   &    69	 &    69	&    69  &0.2648 \\
			SENSORS   &   100	&    26	&    26	&    26	&0.1197   &    49	 &    49	&    49  &0.2018 \\
			SINEVAL   &     2	&   236	&   167	&   167	&0.0095   &   123	 &    89	&    89  &0.0068 \\
			SINQUAD   &  5000	&    85	&    76	&    76	&0.2383   &    80	 &    73	&    73  &0.1903 \\
			SISSER   &     2	&    13	&    13	&    13	&0.0015   &    13	 &    13	&    13  &0.0011 \\
			SNAIL   &     2	&     9	&     9	&     9	&0.0030   &     9	 &     9	&     9  &0.0008 \\
			SPARSINE   &  5000	& 11847	& 10201	& 10201	&29.1683   & 10975	 &  9761	&  9761  &26.5816 \\
			SPARSQUR   & 10000	&    32	&    32	&    32	&0.1762   &    34	 &    34	&    34  &0.1928 \\
			SPMSRTLS   &  4999	&   757	&   667	&   667	&1.5302   &   629	 &   598	&   598  &1.4918 \\
			SROSENBR   &  5000	&    18	&    17	&    17	&0.0302   &    19	 &    18	&    18  &0.0162 \\
			SSBRYBND   &  5000	& 12214	& 10594	& 10594	&29.9288   & 11862	 & 10496	& 10496  &28.9257 \\
			SSI   &     3	& 75345	& 34049	& 34049	&1.6353   &   838	 &   740	&   740  &0.0382 \\
			TESTQUAD   &  5000	&  7117	&  6201	&  6201	&3.5635   &  9752	 &  8542	&  8542  &4.7884 \\
			TOINTGOR   &    50	&   187	&   177	&   177	&0.0207   &   186	 &   186	&   186  &0.0119 \\
			TOINTGSS   &  5000	&     2	&     2	&     2	&0.0065   &     2	 &     2	&     2  &0.0026 \\
			TOINTPSP   &    50	&   223	&   198	&   198	&0.0228   &   204	 &   196	&   196  &0.0098 \\
			TOINTQOR   &    50	&    49	&    49	&    49	&0.0029   &    53	 &    53	&    53  &0.0095 \\
			TQUARTIC   &  5000	&    50	&    27	&    27	&0.0411   &    32	 &    27	&    27  &0.0349 \\
			TRIDIA   &  5000	&  2477	&  2196	&  2196	&1.6461   &  2821	 &  2573	&  2573  &2.0615 \\
			VARDIM   &   200	&  1986	&   205	&   205	&0.0506   &   757	 &   102	&   102  &0.0147 \\
			VAREIGVL   &    50	&    28	&    28	&    28	&0.0024   &    27	 &    27	&    27  &0.0027 \\
			WATSON   &    12	&  3055	&  2641	&  2641	&0.1339   &  2611	 &  2225	&  2225  &0.1570 \\
			WOODS   &  4000	&   280	&   238	&   238	&0.2771   &   312	 &   285	&   285  &0.2680 \\
			YFITU   &     3	&  4675	&  4079	&  4079	&0.2099   &   660	 &   634	&   634  &0.0294 \\
			ZANGWIL2   &     2	&     3	&     2	&     2	&0.0005   &     3	 &     2	&     2  &0.0003 \\

			\hline
			
		\end{tabular}		
	\end{scriptsize}	
\end{table}

In Tables \ref{tbunc}--\ref{tbunc3}, we present the total number of function evaluations ``nef'' and gradient evaluations ``ngrad'', the number of iterations ``iter'', and the CPU time in seconds ``time'' costed by Algorithm \ref{alunc} and the BBQ method, which show that Algorithm \ref{alunc} is faster than the BBQ method on most of the problems.

\section{Conclusion and discussion}\label{secclu}
A mechanism has been introduced for the gradient method to achieve three-dimensional quadratic termination. Then we utilized it to derive a novel stepsize $\alpha_k^{new}$ such that the family \eqref{fystepbb} can have the three-dimensional quadratic termination property by equipping $\alpha_k^{new}$. When applied to the Barzilai--Borwein (BB) method, the computation of $\alpha_k^{new}$ only requires stepsizes and gradient norms in recent three iterations and does not need any exact line search or the Hessian. An efficient gradient method, the method \eqref{snew}, adaptively using short steps associated with $\alpha_k^{new}$ and long BB steps, for unconstrained quadratic optimization is developed. Based on the method \eqref{snew}, we proposed an efficient method, Algorithm \ref{alunc}, for unconstrained optimization. Numerical experiments show evidence of the proposed method over recent successful gradient methods in the literature.

The proposed mechanism can be easily extended to the four-dimensional case. Indeed, let $u_k, v_k, r_k, w_k\in\mathbb{R}^n$ be orthonormal vectors and
\begin{equation*} \label{eqhesgft4}
	H_k
	=\begin{pmatrix}
		u_k\tr  A u_k & u_k\tr A v_k & u_k\tr A r_k & u_k\tr A w_k\\
		v_k\tr A u_k   & v_k\tr A v_k & v_k\tr A r_k & v_k\tr A w_k\\
		r_k\tr A u_k   & r_k\tr A v_k & r_k\tr A r_k & r_k\tr A w_k\\
		w_k\tr A u_k   & w_k\tr A v_k & w_k\tr A r_k & w_k\tr A w_k\\
	\end{pmatrix}.
\end{equation*} 
Similar to \eqref{3dtrdeteqs2}, we have
\begin{equation*}\label{4dtrdeteqs2}
	\begin{aligned}
		\mathrm{tr}(H_k)&=\lambda_1(H_k)+\lambda_2(H_k)+\lambda_3(H_k)+\lambda_4(H_k),\\
		\mathrm{det}(H_k)&=\lambda_1(H_k)\lambda_2(H_k)\lambda_3(H_k)\lambda_4(H_k),\\
		\mathrm{tr}(H_k^2)&=\lambda_1^2(H_k)+\lambda_2^2(H_k)+\lambda_3^2(H_k)+\lambda_4^2(H_k),\\
		\mathrm{tr}(H_k^3)&=\lambda_1^3(H_k)+\lambda_2^3(H_k)+\lambda_3^3(H_k)+\lambda_4^3(H_k).
	\end{aligned}
\end{equation*}
It is not difficult to see that the largest root of
\begin{equation}\label{qt4eq4}
	z^4-\mathrm{tr}(H_k)z^3+\frac{\mathrm{tr}^2(H_k)-\mathrm{tr}(H_k^2)}{2}z^2-\vartheta z+\mathrm{det}(H_k)=0
\end{equation}
is $\lambda_4(H_k)$, where
\begin{equation*}
	\vartheta=\frac{\mathrm{tr}^3(H_k)+2\mathrm{tr}(H_k^3)-3\mathrm{tr}(H_k)\mathrm{tr}(H_k^2)}{6}.
\end{equation*}
Thus, using similar arguments as those in Theorem \ref{3dftth}, the reciprocal of the largest root of \eqref{qt4eq4} together with $\alpha_k^{new}$ and some stepsize that has the two-dimensional quadratic termination property will yield four-dimensional quadratic termination.

Our results further demonstrate potential benefits of the higher-dimensional quadratic optimization model as compared to the two-dimensional one. An interesting question is whether higher-dimensional quadratic optimization models will be helpful in designing efficient gradient methods for constrained optimization.


\bibliographystyle{siamplain}

\end{document}